\documentclass{amsart} 

\usepackage{eucal, mathrsfs}
\usepackage[OT2,OT1,T1]{fontenc}
\usepackage{amsmath,amssymb,amstext,amsthm,amsfonts,amscd,amsopn,epic,eepic}
\usepackage{hyperref}
\usepackage{newcent}
\usepackage{epsfig}
\usepackage{xspace}
\usepackage{epsfig,latexsym,color}
\usepackage[all]{xy}

\theoremstyle{plain}
\newtheorem{thm}{Theorem}[section] 
\newtheorem{pro}[thm]{Proposition}
\newtheorem{lem}[thm]{Lemma}
\newtheorem{cor}[thm]{Corollary}
\newtheorem*{thm2}{Theorem}

\newtheorem*{lem2}{Lemma}

\theoremstyle{definition}
\newtheorem*{defi}{Definition}

\theoremstyle{remark}
\newtheorem*{rem}{Remark}

\newtheorem*{ex}{Example}

\newtheorem{proclaim}{Theorem}
\newtheorem{example}[proclaim]{Example} 

\def\pf{{\em Proof.}\ \,}

\def\({\left(}
\def\){\right)}

\def\nt{\noindent}

\def\fun{\rightarrow}
\def\lfun{{\longrightarrow}}

\def\bC{{\mathbb{C}}}
\def\bD{{\mathbb{D}}}

\def\bP{{\mathbb{P}}}
\def\bQ{{\mathbb{Q}}}
\def\bR{{\mathbb{R}}}
\def\bZ{{\mathbb{Z}}}

\def\cA{{\mathcal A}}

\def\cD{{\mathcal D}}
\def\cE{{\mathcal E}}
\def\cExt{{\mathcal E}xt}
\def\cF{{\mathcal F}}

\def\cH{{\mathcal H}}
\def\cHom{{\mathcal H}om}
\def\cI{{\mathcal I}}
\def\cL{{\mathcal L}}
\def\cM{{\mathcal M}}
 
\def\cO{{\mathcal O}}
\def\cP{{\mathcal P}}
\def\cT{{\mathcal T}}

\def\cU{{\mathcal U}}
\def\cZ{{\mathcal Z}}

\def\ms{\mathscr}

\def\ch{{\textrm{ch}}}
\def\codim{{\textrm{codim\,}}}
\def\coker{{\textrm{coker\,}}}

\def\Ext{{\textrm{Ext\,}}}

\def\Hom{{\textrm{Hom\,}}}

\def\id{{\textrm{id}}}
\def\Im{{\textrm{Im\,}}}
\def\NS{{\textrm{NS}}}
\def\Pic{{\textrm{Pic}}}

\def\Real{{\textrm{Re}}}
\def\rk{{\textrm{rk}}}

\def\Spec{{\textrm{Spec\,}}}

\begin{document}

\title{Bridgeland-stable Moduli Spaces for K-trivial surfaces}


\author{Daniele Arcara}


\address{Department of Mathematics, Saint Vincent College,
300 Fraser Purchase Road, Latrobe, PA 15650-2690, USA}


\email{daniele.arcara@email.stvincent.edu}



\author{Aaron Bertram}

\address{Department of Mathematics, University of Utah,
155 S. 1400 E., Room 233, Salt Lake City, UT 84112-0090, USA}
\email{bertram@math.utah.edu}



\address{Fine Hall, Washington Road, Princeton NJ 08544-1000}
\email{lieblich@math.princeton.edu}

\thanks{The first author was partially supported by a faculty research grant
from St. Vincent College. The second author was partially supported by NSF grant
DMS 0501000}


\begin{abstract}

We give a natural family of Bridgeland stability conditions on the derived
category of a 
smooth projective complex surface $S$ and describe ``wall-crossing behavior''
for objects with the same invariants as $\cO_C(H)$ when 
$H$ generates Pic$(S)$
and $C \in |H|$.
If, in addition,  $S$ is a K3 or Abelian surface, we use this description to construct a 
sequence of fine moduli spaces of Bridgeland-stable objects via Mukai flops and generalized 
elementary modifications of the universal coherent sheaf. We also discover a natural generalization
of Thaddeus' stable pairs for curves embedded in the moduli spaces.

\end{abstract} 

\maketitle

\tableofcontents

\section{Introduction}

Let $X$ be a smooth complex projective variety of dimension $n$.
An ample divisor class $H$ on $X$ defines a slope function on torsion-free sheaves $E$  on $X$ via:
$$\mu_H(E) = \left(\int_X c_1(E) \cdot H^{n-1} \right)/\rk(E)$$
This slope function is a measure of the growth of the Hilbert function of $E$, but it also 
allows one to define the important notion of $H$-stability:

\medskip

\nt {\bf Definition:} $E$ is $H$-stable if $\mu_H(F) < \mu_H(E)$ for all $F \subset E$ with $\rk(F) < \rk(E)$. 

\medskip

\nt It is well-known that this notion allows one to classify the torsion-free sheaves  on $X$ via:

\medskip

$\bullet$ Moduli  for $H$-stable torsion-free sheaves with fixed Hilbert polynomial,

\medskip

$\bullet$ Jordan-H\"older filtrations for $H$-semi-stable torsion-free sheaves, and

\medskip

$\bullet$ Harder-Narsimhan filtrations for arbitrary torsion-free sheaves.

\medskip

A {\it Bridgeland slope function}, in contrast, is defined on the (bounded) derived category 
$\cD(X)$ of coherent sheaves on $X$. It is a pair $(Z,\cA^\#)$ consisting of a linear {\it central charge:}
$$Z:  K(\cD(X)) \rightarrow \bC$$
on the Grothendieck group, together with the heart $\cA^\#$ of a t-structure on 
$\cD(X)$ that is compatible with the central charge in the sense that
$$Z(A)  \in \cH = \{ \rho e^{i\phi} \ | \rho > 0, \ 0 < \phi \le \pi\}$$
for all non-zero objects $A$ of $\cA^\#$. This allows one to define a (possibly infinite-valued) slope:
$$\mu_Z(A) := -\frac{\mbox{Re}(Z(A))}{\Im(Z(A))}$$
for objects of $\cA^\#$ analogous to the $H$-slope on coherent sheaves. The pair 
$(Z,\cA^\#)$ is called a {\it Bridgeland stability condition} if the associated notion of $Z$-stability 
has the Harder-Narasimhan property.

\medskip

In this paper, we will consider central charges of the form:
$$Z([E]) = -\int_S e^{-(D + iF)} \ch([E]) \ \ \mbox{and}\ \ \ Z'([E]) = -\int_S e^{-(D + iF)} \ch([E]) \sqrt {\mbox{td}(S)}$$
on a smooth projective surface $S$, where $F$ is an ample $\bR$-divisor, and $D$ is an arbitrary 
$\bR$-divisor. Following Bridgeland's argument for $K3$ surfaces \cite{Bri03}, we show that the former
always has a natural partner $t$-structure $\cA^\#$ (depending upon $D$ and the {\bf ray} generated by $F$) such that the pair $(Z,\cA^\#)$ defines a stability condition.
Our main results focus further on the one-parameter family of stability conditions on a fixed abelian category 
$\cA^\#$, where:
$$\mbox{Pic}(S) = \bZ[H], \ \ D = \frac 12 H, \ \mbox{and}\ \ F = tH; \ \ t > 0$$

This family of stability conditions is well-tuned to study the stability of objects $E \in \cA^\#$ with chern class invariants:
$$\ch(E) = H + H^2/2 = \ch(\cO_S(H)) - \ch(\cO_S)$$
in the sense that we will be able to state the precise set of stable objects (depending on $t$) with those invariants.
Moreover, in the $K$-trivial case (i.e. when $S$ is a K3 or Abelian surface), we will use this knowledge to 
construct proper moduli spaces of Bridgeland-stable objects by starting with the relative Jacobian (the moduli of 
stable objects for $t >> 0$) and performing a sequence of Mukai flops as $t$ passes over a series
of ``walls.'' This in particular exhibits a sequence of birational models of the relative Jacobian, 
which seem to be new, although they encode quite a lot of interesting results on the positivity of the line bundle 
$\cO_S(H)$ on the surface.

\medskip

To get an idea of the wall-crossing phenomenon, consider the exact sequence of coherent sheaves:
$$0 \fun \cO_S \fun \cO_S(H) \fun i_*\cO_C(H)  \fun 0$$ 
for some curve $C \in |H|$. This, it turns out, will {\bf not} be an exact sequence of objects 
in our category $\cA^\#$. Rather, 
the sequence:
$$0 \fun \cO_S(H) \fun  i_*\cO_C(H) \fun \cO_S[1] \fun 0$$
coming from the ``turned'' distinguished triangle in $\cD(S)$ will be a short exact sequence of objects of $\cA^\#$. 
Below the critical ``wall'' value $t = \frac 12$, we will have $\mu_Z(\cO_S(H)) > \mu_Z(i_*\cO_C(H))$, exhibiting 
$i_*\cO_C(H)$ as an unstable object of $\cA^\#$(!). The ``replacement'' stable object(s) will be of the form:
$$0 \fun \cO_S[1] \fun  E \fun \cO_S(H) \fun 0$$
which are parametrized by $\bP(\Ext^1_{\cA^\#}(\cO_S(H),\cO_S[1])) = \bP(\Ext^2_S(\cO_S(H),\cO_S)) \cong
\bP(H^0(S,\cO_S(H))^*)$ via Serre duality. 

\medskip

Our moduli functor is based upon the  generalized notion of a flat family we learned from
Abramovich and Polishchuk \cite{AP06}. One would, of course, like to have an a priori construction of 
moduli spaces of 
Bridgeland-stable objects via some sort of invariant theory argument, but 
the fact that we are not working exclusively with coherent sheaves makes it difficult to
see how to make such a construction. Instead, we rely on the fact that an Artin stack of  flat families of
objects of $\cA^\#$ exists, using a result of Max Lieblich, which we attach as an appendix, 
and then work rather hard to show that stability is an open condition in the cases of interest to us
(recent work by Toda \cite{Tod07} gives an alternative, and more general, approach).
We then work by induction, starting with the universal family over the relative Jacobian and 
elementary modifications across the Mukai flops to actually prove that each successive birationally
equivalent space is indeed a fine moduli space of Bridgeland stable objects. Thus we are able in this 
case to carry out the program envisioned by Bridgeland at the very end of \cite{Bri03}.

\medskip

The methods introduced here should be useful in the construction of Bridgeland stable 
moduli spaces of objects with 
other invariants on surfaces both with and without the $K$-trivial assumption. In particular, in
joint work with Gueorgui Todorov \cite{ABT}, we will describe the menagerie of Mukai flops of Hilbert schemes of 
$K$-trivial surfaces induced by varying Bridgeland stability conditions on objects with 
the invariants of an ideal sheaf of points.

\medskip

\nt {\bf Acknowledgements.} We thank 
Dan Abramovich, Tom Bridgeland, Andrei C\u ald\u araru, Max Lieblich, Dragan Mili\v ci\' c
and Alexander Polishchuk
for all of their help, especially on the subjects of Bridgeland stability conditions and the
subtleties of the derived
category. 

\medskip

\section{Stability Conditions on the Derived Category of a Surface} 

We start with some general remarks on the bounded derived category of 
coherent sheaves $\cD(S)$ for the uninitiated reader.
Derived categories were introduced by Verdier in 
\cite{Ver63}. For a comprehensive introduction, see \cite{Mil}.

\medskip

The objects of $\cD(S)$ are complexes (with bounded cohomology):
$$ \cdots \lfun E_{i-1} \overset{d_{i-1}}{\lfun} E_i
\overset{d_i}{\lfun} E_{i+1} \lfun \cdots $$
of coherent sheaves, and homotopy classes of maps of complexes are maps in $\cD(S)$. Let:
$$ \cH^i(E) := \ker d_i / \textrm{im} \, d_{i-1} $$
denote the cohomology sheaves of the complex. A (homotopy class of) map(s) of complexes:
\begin{eqnarray*}
\begin{CD}
\cdots @>>> E_{i-1} @> d_{i-1} >> E_i @> d_i >> E_{i+1} @>>> \cdots \\
& & @V f_{i-1} VV @V f_i VV @V f_{i+1} VV \\
\cdots @>>> F_{i-1} @> d_{i-1} >> F_i @> d_i >> F_{i+1} @>>> \cdots
\end{CD}
\end{eqnarray*}
is a quasi-isomorphism if it induces isomorphisms on all cohomology sheaves. The full 
class of maps in $\cD(S)$ is obtained by formally inverting all quasi-isomorphisms.
Thus, in particular, quasi-isomorphic complexes represent isomorphic objects.

\medskip

Here are a few \medskip facts about the derived category:

\medskip

$\bullet$ $E[1]$ (and $f[1]$) denote the shifts: $(E[1])_i = E_{i+1}, \  (f[1])_i = f_{i+1}. $

\medskip

$\bullet$ $\cD(S)$ is a triangulated category. It does not make sense to talk about kernels and
cokernels of a map $f:E \rightarrow F$. Rather, the map $f$ induces a 
cone $C$ and a {\it distinguished
triangle}
$$ \cdots \lfun E \overset{f}{\lfun} F \lfun C \lfun E[1]
\overset{f[1]}{\lfun} F[1] \lfun \cdots $$

\medskip

$\bullet$ Given two coherent sheaves $E$ and $F$, there is an isomorphism
$$ \Hom_{\cD(S)}(E[m], F[n]) \cong \Ext^{n-m}_S(E, F), $$

\medskip

$\bullet$ A short exact sequence of sheaves
$$0 \lfun K \lfun E \lfun Q \lfun 0$$
induces a distinguished triangle
$$ \cdots \lfun K \lfun E \lfun Q \lfun K[1] \lfun E[1] \fun Q[1] \fun \cdots$$
where  $(Q\fun K[1])\in\Hom_{\cD(S)}(Q,K[1])=\Ext^1_S(Q,K)$ is 
the extension class.

\medskip

$\bullet$ A distinguished triangle
$$ \cdots \lfun F \lfun E \lfun G \lfun F[1] \lfun \cdots $$
induces a long exact sequence of cohomologies
$$ \cdots \lfun \cH^i(F) \lfun \cH^i(E) \lfun \cH^i(G) \lfun \cH^{i+1}(F)
\lfun \cdots $$
in the category of coherent sheaves on $S$.

\medskip

$\bullet$ The derived dual $E^\vee$ of an object $E\in\cD(S)$ is defined as
$R\cHom(E,\cO_S)$, where $R\cHom$ is the derived functor induced by the
$\cHom$ functor on coherent sheaves.

\medskip

$\bullet$
It is always true that $E^\vee{}^\vee=E$ for objects of the derived category.

\medskip

$\bullet$ If $E$ is a sheaf, the derived dual $E^\vee$ is represented by a complex with
$$ \cH^i(E^\vee) = \cExt^i_S(E,\cO_S). $$

\medskip

$\bullet$ The derived dual of a distinguished triangle
$$ \cdots \lfun F \lfun E \lfun G \lfun F[1] \lfun \cdots $$
is the distinguished triangle
$$\cdots \lfun G^\vee \lfun E^\vee \lfun F^\vee \lfun G^\vee[1] \lfun \cdots$$

\medskip

$\bullet$ If three consecutive terms $F$, $E$,
and $G$ of a distinguished triangle 
$$ \cdots \lfun F \lfun E \lfun G \lfun F[1] \lfun \cdots $$
are in the heart, $\cA$, of a $t$-structure on $\cD(S)$, then they determine a short exact sequence
$0 \lfun F \lfun E \lfun G \lfun 0$ of objects in $\cA$.

We will only be concerned here with particular sorts of $t$-structures on $\cD(S)$ obtained by 
{\it tilting}. In general, tilting is obtained as follows, starting with an abelian category $\cA$.

\begin{defi}[{[HRS96]}] A pair of full subcategories $(\cT,\cF)$ in $\cA$
is called a {\it torsion pair} if:

\medskip

(TP1) $\Hom_{\cA}(T,F)=0$ for every $T\in \cT$ and $F\in \cF$.

\medskip

(TP2) Every object $E\in \cA$ fits into a short exact sequence
$$ 0 \fun T \fun E \fun F \fun 0 \ \ \mbox{with}\ \ T\in \cT \ \ \mbox{and}\ \ F\in \cF$$
\end{defi}

\begin{lem2}[{[HRS96]}]

Let $(\cT,\cF)$ be a torsion pair in an abelian category $\cA$.
If $\cA$ is the heart of a bounded $t$-structure on a triangulated category
$\cD$, then the full subcategory of $\cD$ with objects:
$$ob( \cA^\#) = \{ E \in \cD \mid 
H^{-1}(E) \in \cF, \ H^0(E) \in \cT, H^j(E) = 0 \textrm{ for } j \neq -1, 0 \}, $$
is the heart of another $t$-structure on $\cD$, hence in particular an abelian category.

\end{lem2}

\begin{rem} A short exact sequence:
$$0\rightarrow K \rightarrow E \rightarrow Q \rightarrow 0$$
of objects of $\cA^\#$ gives rise to a long exact sequence of objects of $\cA$:
$$0 \rightarrow H^{-1}(K) \rightarrow H^{-1}(E) \rightarrow H^{-1}(Q) \rightarrow H^0(K) 
\rightarrow H^0(E) \rightarrow H^0(Q) \rightarrow 0$$
where $H^{-1}(K), H^{-1}(E),H^{-1}(Q) \in \cF$ and $H^{0}(K), H^{0}(E),H^{0}(Q) \in \cT$. 
\end{rem}

\noindent {\bf Our Tilts:} Given $\bR$-divisors $D,F \in H^{1,1}(S,\bR)$ on a surface $S$, with $F$ ample, 
then as in \S1, the $F$-slope of a torsion-free coherent sheaf $E$ on $S$ 
is given by:
$$\mu_F(E) = \left({\int_S c_1(E)\cdot F}\right)/{\mbox{rk}(E)}$$
and all coherent sheaves on $S$ have a unique {\it Harder-Narasimhan
filtration}:
$$E_0 \subset E_1 \subset \cdots \subset E_{n(E)} = E$$
characterized by the property that $E_0 = \mbox{tors}(E)$ and each $E_{i}/E_{i-1}$ is a torsion-free 
$F$-semistable sheaf of 
slope $\mu_i$ (i.e. an extension of $F$-stable sheaves of slope $\mu_i$) for a strictly decreasing
sequence
$$\mu_{F-max}(E) := \mu_1 > \mu_2 > \cdots > \mu_{n(E)} =: \mu_{F-min}(E)$$

\begin{defi} Let $\cA$ be the category of coherent sheaves on $S$, and let: 
$$ob(\cT) = \{ \mbox{torsion sheaves}\} \cup \left\{ E \ | \
\mu_{F-min}(E) > \int_S D.F\right\} \ \mbox{and}$$
$$ob(\cF) = \left\{ E \ | \ 
\mu_{F-max}(E) \le \int_S D.F \right\}$$
Note that this only depends upon the ray spanned by $F$. Now 
define $\cA^\#_{(D,F)}$ by applying the Lemma to  the standard $t$-structure on the bounded
derived category of coherent sheaves on $S$.  

\end{defi}

\noindent {\bf Our Charges:} For now we will only consider the charges:
$$Z_{(D,F)}(E) := -\int_S e^{-(D+iF)}\ch(E)$$
extended to $K(\cD)$ by defining:
$Z_{(D,F)}(E) = \sum (-1)^iZ_{(D,F)}(H^i(E))$
for all objects $E$ of $\cD$, so that in particular:
$$Z_{(D,F)}(E) = Z_{(D,F)}(H^0(E)) - Z_{(D,F)}(H^{-1}(E))$$
for objects $E$ of the category $\cA^\#$.

\medskip

Recall the Hodge Index Theorem and the Bogomolov-Gieseker Inequality for surfaces
(see, for example, \cite{Fri98}):

\medskip

\begin{thm2}[Hodge Index] \label{Hodge Index Theorem}
If $D$ is an $\bR$-divisor on $S$ and $F$ is an ample $\bR$-divisor, then:
$$D\cdot F = 0 \Rightarrow D^2 \le 0$$
\end{thm2}

\begin{thm2}[Bogomolov-Gieseker Inequality]\label{Bogomolov-Gieseker Inequality}
If $E$ is an $F$-stable torsion-free sheaf on $S$, then:
$$\emph{ch}_2(E) \le \frac{c_1^2(E)}{2\cdot \emph{rk}(E)}$$
\end{thm2}

As an immediate corollary of these two results, we have:

\medskip

\begin{cor}\label{slope} Each pair $(Z_{(D,F)},\cA^\#_{(D,F)})$ is a Bridgeland slope function.
\end{cor}

\pf Since each object $E$ of $\cA^\#_{(D,F)}$ fits into an exact sequence:
$$0 \fun H^{-1}(E)[1] \fun E \fun H^0(E) \fun 0$$
with $H^{-1}(E) \in \cF$ and $H^0(E) \in \cT$, and since $\cH$ is closed under addition, it suffices to 
show that:

\medskip

(1) $Z_{(D,F)}(T) \in \cH$ for all torsion sheaves on $S$,

\medskip

(2) $Z_{(D,F)}(E) \in \cH$ for all $F$-stable torsion-free sheaves with $\mu_F(E) > D\cdot F$

\medskip

(3) $Z_{(D,F)}(E[1]) \in \cH$ for all $F$-stable torsion-free sheaves with $\mu_F(E) \le D\cdot F$.

\medskip

Let  $Z(E) := Z_{(D,F)}(E)$ and compute:
$$Z(E)  = \left( -\ch_2(E) + D\cdot c_1(E) - \rk(E)\left(D^2/2 - F^2/2\right) \right)+ iF\cdot\left(c_1(E) - \rk(E)D\right)$$

In (1), either $T$ is supported in dimension $1$ and 
$\Im(Z(T)) = c_1(E) \cdot F > 0$ since $c_1(E)$ is
effective, or else $T$ is supported in dimension $0$, in which case $\Im(Z(T)) = 0$, but 
Re$(Z(T)) = - \ch_2(T) < 0$.  So $Z(T) \in \cH$.

\medskip

In (2), $\Im(Z(E)) = F\cdot (c_1(E) - \rk(E)D) = \rk(E) \left(\mu_F(E) - D\cdot F\right) > 0$.
Similarly, in (3), if $\mu_F(E) < D\cdot H$, then $\Im(Z(E)) < 0$, so $\Im Z(E[1]) > 0$. Finally,
if $\mu_F(E) = D\cdot F$ and $E$ is $F$-stable, then by the Bogomolov-Gieseker inequality:
$$\mbox{Re}(Z(E)) \ge -\frac{{c_1}^2(E)}{2\cdot \rk(E)} + D\cdot c_1(E) - \rk(E)\left(D^2/2 - F^2/2\right) $$
$$  \ \ \ \ \ \ \ \ \ \ = - \frac{1}{2\cdot \rk(E)}\left( \rk(E)D - c_1(E)\right)^2 + \rk(E) F^2/2$$
But $\mu_F(E) = D\cdot F$ implies $(\rk(E)D - c_1(E))\cdot F = 0$, and so by the Hodge index theorem,
we have $\mbox{Re}(Z(E)) > 0$, $\mbox{Re}(Z(E[1])) < 0$, and $Z(E[1]) \in \cH$, as desired. \qed

\medskip

\nt {\it Remark.} One can show that in fact each $(Z_{(D,F)},\cA^\#_{(D,F)})$ is a 
Bridgeland stability condition. The argument is the same as in Bridgeland's K3 paper, where 
Harder-Narasimhan filtrations are directly  produced when $D$ and $F$ are $\bQ$-divisors 
(Proposition 7.1),
and the general case is deduced by continuity and the structure of the space of stability conditions.
The other  ``standard'' stability properties are easier to see. For example, for  
$(D,F)$-stable objects $A,B \in \cA^\#$, the implication:
$$\mu_Z(A) > \mu _Z(B) \Rightarrow \mbox{Hom}_{\cA^\#}(A,B) = 0$$
is immediate for Bridgeland slope functions, as is the implication:
$\mu_Z(A) = \mu_Z(B) \Rightarrow A \cong B$. Moreover, Hom$_{\cA^\#}(A,A) = \bC \cdot \id$
(Schur's Lemma) also follows easily by considering the induced map of an isomorphism $f$ (and 
$f - \lambda \cdot \id$) on cohomology sheaves.

\medskip

\section{Some Bridgeland-Stable Objects}

It is tricky, in general, to determine which objects of $\cD(S)$ are stable for a arbitrary pairs $(D,F)$.
Fortunately, to determine ``wall-crossing'' phenomena, it is enough to consider a one-parameter 
family $Z_t := Z_{(D_t,F_t)}$ of central charges. This is what we will do, in the special case:
$$\Pic(S) = \bZ[H]$$
for stability conditions:
$$D = \frac 12H, \ \ F = tH; \ \ t > 0$$
so that if we let $Z_t(E) = Z_{(\frac 12 H, tH)}(E)$, then:
$$Z_t(E)  = -\ch_2(E) + \frac 12 c_1(E)\cdot H + \frac{\rk(E)H^2}2\left(t^2 - \frac 14\right) + 
i t\left(H\cdot \left(c_1(E) - \frac {\rk(E)H}2\right)\right)$$

The abelian category $\cA^\# := \cA^\#_{(\frac 12 H,tH)}$ is independent of $t$, and the 
$t$-stable objects of $\cA^\#$ of infinite $Z_t$-slope are always either of the form:

\medskip

(a) Torsion sheaves supported in dimension zero, or

\medskip

(b) Shifts $E[1]$ of $H$-stable vector bundles $E$ of even rank $2n$ and $c_1(E) = nH$.

\medskip

\noindent All objects in $\cA^\#$ of infinite slope are extensions of these. For example, if $F$ is 
torsion-free but not locally free of rank $2n$ and 
$c_1(F) = nH$, let $E = F^{**}$ and then:
$$0 \rightarrow F \rightarrow E \rightarrow T_Z \rightarrow 0$$
(for a sheaf $T_Z$ supported on a scheme $Z$ of dimension zero) becomes the exact sequence:
$$0 \rightarrow T_Z \rightarrow F[1] \rightarrow E[1] \rightarrow 0$$
exhibiting $F[1]$ as an extension in $\cA^\#$ of $E[1]$ by the torsion sheaf $T_Z$.

\begin{defi} An object $E$ of $\cA^\#$ is {\it $t$-stable} if it is stable with respect to the
central charge
$Z_t$, i.e.
$$\mu_t(K) = - \frac{{\mbox{Re}}(Z_t(K))}{\Im(Z_t(K))}< \mu_t(E)$$
for all subobjects $K \subset E$ in $\cA^\#$ ($\Leftrightarrow \mu_t (E) < \mu_t(Q)$
for all surjections $E \fun Q \fun 0$ in $\cA^\#$).
\end{defi}

The following Lemma establishes the $t$-stability of some basic objects of
$\cA^\#$.

\medskip

\begin{lem}\label{LambdaI} Let $Z \subset S$ be a subscheme of dimension zero. Then the objects:
$${\cI}_Z(H) \  \ \ \ \ \mbox{and}\ \ \ \ \ \ \cI_Z^\vee[1] = RHom(\cI_Z,\cO_S)[1]$$
are $t$-stable for all $t > 0$.
\end{lem}

\pf  Let $K$ be a subobject of $\cI_Z(H)$ and let
$$(*) \ \  0 \fun K \fun \cI_Z(H) \fun Q \fun 0 $$
be the associated  ``potentially destabilizing'' sequence in $\cA^\#$.
This induces:
$$0 \fun H^{-1}(K) \fun 0 \fun H^{-1}(Q) \fun H^0(K) \fun \cI_Z(H) \fun H^0(Q) \fun 0$$
which imples that $H^0(Q)=\cI_Z(H)$ or else $H^0(Q)$ is a torsion sheaf.

\medskip

If $H^0(Q)=\cI_Z(H)$, then $H^0(K) = H^{-1}(Q) \in \cT \cap \cF = 0$, so $K = 0$.

\medskip
 
If $H^0(Q)$ is torsion,
then $H^0(K) \in \cT$ and $H^{-1}(Q) \in \cF$ imply that the support of 
$H^0(Q)$ has dimension $0$,  and $\rk(H^0(K)) = \rk(H^{-1}(Q)) + 1$ and $c_1(H^0(K)) = c_1(H^{-1}(Q)) + H$, so
$\rk(H^{-1}(Q)) = 2n$ and $c_1(H^{-1}(Q)) = nH$ (or else $H^{-1}(Q) = 0$). Thus
$Q$ has infinite slope(!), and the ``potentially destabilizing'' sequence
cannot, therefore, destabilize $\cI_Z(H)$.

\medskip

Turning next to $\cI_Z^\vee[1]$, notice that:
$$H^{-1}(\cI_Z^\vee[1]) = \cO_S \ \mbox{and}\ \ H^0(\cI_Z^\vee[1]) = \cExt^2_{\cO_S}(\cO_Z,\cO_S) = T
\ \mbox{(torsion, supported on $Z$)}$$
so any ``potentially destabilizing'' sequence:
$$(*) \ \ 0 \rightarrow K \rightarrow \cI_Z^\vee[1] \rightarrow Q \rightarrow 0$$
gives rise to a long exact sequence of coherent sheaves:
$$ 0 \fun H^{-1}(K) \fun \cO_S \fun H^{-1}(Q) \fun H^0(K) \fun T \fun H^0(Q) \fun 0$$

Thus $H^0(Q)$ is supported in dimension $0$, and  either $H^{-1}(K) = 0$
or else $H^{-1}(K) = \cO_S$ (otherwise the sheaf $H^{-1}(Q)$ would have
torsion, which is not allowed). 

 If $H^{-1}(K) = \cO_S$, then $c_1(H^{-1}(Q)) = c_1(H^0(K))$ and $\rk(H^{-1}(Q)) = \rk(H^0(K))$, which
contradicts $H^{-1}(Q) \in \cF$ and $H^0(K) \in \cT$, unless of course
$H^{-1}(Q) = 0$. This would not destabilize $\cI_Z^\vee[1]$, since 
$Q = H^0(Q)$ would have infinite slope, since it would be a torsion sheaf
supported in dimension zero.
 
If $H^{-1}(K) = 0$, then either $H^{-1}(Q)$ is locally free of rank $2n$ and $c_1 =
nH$, or else $H^{-1}(Q) = \cO_S$ and $H^0(K)$ is torsion, supported in dimension $0$. 
But in the first case,
$H^{-1}(Q)[1]$ and
$H^0(Q)$ have infinite slope, so
$Q$ has infinite slope and $(*)$ does not destabilize. In the second case, we
need to worry about $H^0(K)$. If
$H^0(K) \ne 0$, then $(*)$
would destabilize $\cI_Z^\vee[1]$ because $K = H^0(K) \subset \cI_Z^\vee[1]$ would have
infinite slope. But the derived dual
contains no such sub-objects $K$, because $(\cI_Z^\vee)^\vee = \cI_Z$ is a sheaf!  \qed

 \begin{rem} The object $\cI_Z^\vee[1]$ is a surface analogue of the line bundle
$\cO_C(D) = I_D^\vee$ on a curve:
$$0 \rightarrow \cO_C \rightarrow \cO_C(D) \rightarrow \cExt^1_{\cO_C}(\cO_D,\cO_C) \rightarrow 0$$
\end{rem}

\begin{rem} Beware of the temptation to treat $t$-stability too casually! Observe:
\end{rem}

\begin{lem}\label{LambdaInot} The objects $\cI_Z[1] \in \cA^\#$ are
not $t$-stable for {\bf any} value of $t$.
\end{lem}

\pf
The exact sequence of objects in $\cA^\#$:
$$0 \rightarrow \cO_Z \rightarrow \cI_Z[1] \rightarrow \cO_S[1] \rightarrow 0$$
always destabilizes $\cI_Z[1]$, since $\cO_Z$ has infinite slope. \qed

\medskip

Next, we turn our attention to the objects of $\cA^\#$ with the invariants:
$$\ch(E) = 0 + H + H^2/2$$
A quick computation shows that:
$$\mu_t(E) = \mbox{Re}(Z_t(E)) = 0 \ \mbox{for all $t$}$$
for all $E$ with these invariants. Thus to check that $\mu_t(K) < \mu_t(E)$ (or $\mu_t(K) > \mu_t(E)$)
it suffices to compute the real part of $Z_t(K)$.

\medskip

\begin{pro}\label{Stable} If $E \in \cA^\#$ has (formal) invariants
$r=0,c_1=H,\emph{ch}_2 = H^2/2$
and $E$ is $t$-stable for {\bf some} value of $t$,
then either $E = i_*L_C$ for some torsion-free rank-one sheaf $L_C$ supported on a curve
$C \in |H|$, or else:

\medskip

(1) $H^0(E)$ has torsion (if any) only in dimension $0$

\medskip

(2) $H^{-1}(E)$ is locally free and $H$-stable of odd rank $2n+1$ with $c_1 = nH$

\medskip

(3) $H^0(E)/tors(H^0(E))$ is $H$-stable of rank $2n+1$ with $c_1 = (n+1)H$.

\medskip

(4) The kernel of $E \rightarrow H^0(E)/tors(H^0(E))$ is the (shifted) derived dual of a torsion-free sheaf.

\end{pro}

\pf If $E$ is a sheaf with these invariants, then it is of the form $i_*L_C$ and any torsion in $L_C$
would destabilize, as all torsion sheaves on $S$ supported in dimension $0$ have infinite slope.
Otherwise
$H^{-1}(E) \ne 0$. Now $\rk(H^{-1}(E)) = \rk(H^0(E))$
and $c_1(H^{-1}(E)) = c_1(H^0(E)) - H$ from the invariants, and this, together with 
$H^{-1}(E) \in \cF$ and $H^0(E) \in \cT$, forces (1). If $\rk(H^{-1}(E)) = 2n$, then 
$c_1(H^{-1}(E)) = nH$ is also forced, and $H^{-1}(E)[1] \subset E$
would be a sub-object of infinite slope, contradicting the stability of $E$ for each value of $t$.
So
$H^{-1}(E)$ has odd  rank $2n+1$, and $c_1(H^{-1}(E)) = nH$. Similarly, $H^{-1}(E)$ and  
$H^0(E)/tors(H^0(E))$ are $H$-stable, and if
$H^{-1}(E)$ were not locally free, then (as in Lemma \ref{LambdaInot}) there would
be a subobject:
$$H^{-1}(E)^{**}/H^{-1}(E) \subset H^{-1}(E)[1] \subset E$$
of infinite slope, contradicting $t$-stability of $E$ for each value of $t$. This gives (2)
and (3). 

\medskip

Finally, let $E'$ be the kernel of the short exact sequence in $\cA^\#$:
$$0 \fun E' \fun E \fun H^0(E)/tors(H^0(E)) \fun 0$$
so:
$H^{-1}(E') = H^{-1}(E)$ is locally free, and $H^0(E') = tors(H^0(E))$ is supported in dimension $0$. 
Then (4) follows directly from:

\begin{lem}\label{Dual}
Suppose $E$ is an object of $\cA^\#$ satisfiying:
$$H^{-1}(E) \ \mbox{is locally free, and}\ H^0(E) \ \mbox{is torsion, supported in dimension}\ 0$$
Then either $E^\vee[1]$ is a torsion-free sheaf or $E$ has a torsion subsheaf supported in dimension $0$.
\end{lem}

\pf Basically, this consists of taking duals twice. First, the dual of the sequence in $\cA^\#$:
$$0 \fun H^{-1}(E)[1] \fun E \fun H^0(E) \fun 0$$
is a distinguished triangle:
$$\cdots \fun H^0(E)^\vee \fun E^\vee \fun H^{-1}(E)^\vee[-1] \fun H^0(E)^\vee[1] \fun \cdots$$
whose associated sequence of cohomology sheaves is:
$$0 \fun H^1(E^\vee) \fun H^{-1}(E)^* \fun H^2(H^0(E)^\vee) \fun H^2(E^\vee) \fun 0$$
because of the assumptions on $H^{-1}(E)$ and $H^0(E)$. Thus either $H^2(E^\vee) = 0$ and $E^\vee[1] = H^1(E^\vee)$ is 
a subsheaf of the dual vector bundle $H^{-1}(E)^*$, or else $H^2(E^\vee) \ne 0$ is the quotient of a sheaf 
supported in dimension $0$, hence is itself a sheaf supported in dimension $0$. But in the latter case,
the distinguished triangle:
$$\cdots \fun H^1(E^\vee)[-1] \fun E^\vee \fun H^2(E^\vee)[-2] \fun H^1(E^\vee) \fun \cdots$$
(coming from the fact that $E^\vee$ has cohomology only in degrees $1$ and $2$) dualizes to:
$$\cdots \fun H^2(E^\vee)^\vee[2] \fun E \fun H^1(E^\vee)^\vee[1] \fun \cdots$$
which is an exact sequence in $\cA^\#$! And if $H^2(E^\vee) \ne 0$, then  $H^2(H^2(E^\vee)^\vee) \ne 0$ as well. 
\qed

\rem Figuring out which of the objects of Proposition \ref{Stable} are $t$-stable for 
each {\bf particular} $t$ is, of course, more delicate. For example:

\begin{lem}\label{StableL} A sheaf $i_*L_C$ as in Proposition \ref{Stable} is
$t$-stable unless there is an $H$-stable torsion-free sheaf $K$ on $S$ of
rank $2n+1$ with $c_1(K) = (n+1)H$, Re$(Z_t(K)) \le 0$ and a map:
$$f: K \rightarrow i_*L_C$$
that is generically (on $C$) surjective with a locally free
kernel.

\end{lem}

\pf On the one hand, such a sheaf (and map) does destabilize $i_*L_C$, since the sequence:
$0 \rightarrow \ker(f) \rightarrow K \rightarrow i_*L_C \rightarrow \cO_Z \rightarrow 0$
is a short exact sequence of objects in $\cA^\#$:
$$0 \rightarrow K \rightarrow i_*L_C \rightarrow Q \rightarrow 0 $$
(because $K \in \cT$ and $F = \ker(f) \in \cF$) where $Q \in \cA$ satisfies
$H^{-1}(Q) = \ker(f)$ and $H^0(Q) = \cO_Z$.
On the other hand, any potentially destabilizing sequence $0 \rightarrow K \rightarrow i_*L_C 
\rightarrow Q \fun 0$ gives rise to:
$$0 \rightarrow H^{-1}(Q) \rightarrow H^0(K) \rightarrow i_*L_C \rightarrow H^0(Q) \rightarrow 0$$

From this we may read off:

\medskip

$\bullet$ $H^0(Q)$ is supported in dimension $0$ (otherwise $H^{-1}(Q) = 0 = H^0(K)$)

\medskip

$\bullet$  $H^0(K)$ is torsion-free

\medskip

$\bullet$ $\rk(H^{-1}(Q)) = r = \rk(H^0(K))$ and $c_1(H^{-1}(Q)) + H = c_1(H^0(K))$. 

\medskip

If $r = 2n$,  then $c_1(H^{-1}(Q)) = nH$, and $Q$ has infinite slope. Thus, if the
sequence is to destabilize $i_*L_C$, it must be the case that $\rk(H^0(K)) = 2n+1$ and $c_1(H^0(K)) =
(n+1)H$, and then it follows from $H^0(K) \in \cT$ that $H^0(K)$ is $H$-stable.  Finally, if
$H^{-1}(Q)$ is  not locally free, then $H^{-1}(Q)^{**}$ has smaller slope, and $Q$ can be replaced by
$Q'$, with:
$$0 \rightarrow H^{-1}(Q)^{**}/H^{-1}(Q) \rightarrow Q \rightarrow Q' \rightarrow 0 \ \mbox{and}\ 
H^{-1}(Q') = H^{-1}(Q)^{**}$$
\qed

\begin{ex} Let $L_C = \cO_C(H + D - D')$ where $D,D'$ are 
effective disjoint divisors of the same degree supported on the smooth part of $C$. 
Then:
$$(*) \ 0 \rightarrow \cO_S \rightarrow \cI_{D'}(H) \rightarrow i_*L_C \rightarrow 
L_C|_D \rightarrow 0$$
(thinking of $D,D'$ as zero-dimensional subschemes of $S$) will $t$-destabilize $i_*L_C$ if:
$$\mbox{Re}(Z_t(\cI_{D'}(H))) =  \deg(D') + t^2\frac{H^2}2 - \frac{H^2}8 \le 0$$ 
or in other words, if:
$$t^2 \le \frac 14 - \deg(D')\frac 2{H^2}$$  
and if equality holds, then $i_*L_C$ will be an extension of stable objects of $\cA^\#$ 
of the same phase (i.e. $i_*L_C$ is $t$-semistable, and $t$ is a critical 
value).

\end{ex} 

One final lemma on $t$-stability will be useful for us:

\medskip

\begin{lem}\label{StableExt} Suppose $E$ is an object of $\cA^\#$ which is an extension (in $\cA^\#$) of the form:
$$0 \fun \cI_W^\vee[1] \fun E \fun \cI_Z(H) \fun 0$$
for zero-dimensional schemes $Z,W$ with $len(W) = len(Z)$ (the rank-one case of
Proposition
\ref{Stable}). Then 
$E$ is $t$-stable unless either Re$(Z_t(\cI_Z(H))) \ge 0$ and the 
quotient 
$E \rightarrow \cI_Z(H)$ destabilizes $E$, or else there is
a sheaf $K \subset E$ as in Lemma \ref{StableL}.
\end{lem}

\pf A potentially $t$-destabilizing sequence $0\rightarrow K \rightarrow E \rightarrow Q
\rightarrow 0$ gives rise to:
$$0 \rightarrow H^{-1}(K) \rightarrow \cO_S \rightarrow H^{-1}(Q) \rightarrow 
H^0(K) \rightarrow H^0(E) \rightarrow H^0(Q) \rightarrow 0$$
The (by now) standard analysis allows us to conclude that if the sequence actually destabilizes $E$,
then one
of the following two must be true. Either:

\medskip

$\bullet$ $H^{-1}(K) = 0$, $H^0(Q)$ is torsion supported in dimension $0$, and $H^0(K) = K$ of the
Lemma, or:

\medskip

$\bullet$ $H^{-1}(K) = 0$, $H^0(Q) = \cI_Z(H)$, and rk$(H^{-1}(Q)) = 2n, c_1(H^{-1}(Q)) = nH$.

\medskip

\noindent But in the latter case, the slope of $H^0(Q) = \cI_Z(H)$ is smaller than the slope of $Q$, so
if $Q$ destabilizes $E$, then so does $\cI_Z(H)$ (only more so!).

\qed

\medskip

\begin{thm}\label{tstable}The objects $E$ of $\cA^\#$ with numerical invariants:
$$\emph{ch}_0(E) = 0,\ \emph{ch}_1(E) = H, \ 
\emph{ch}_2(E) = \frac{H^2}2$$
that are $t$-stable for
some $t > \frac 16$ are either:

\medskip

$\bullet$ Sheaves of the form $i_*L_C$ (as in Proposition \ref{Stable}), or else

\medskip

$\bullet$ Fit into (non-split!) extensions of the form:
$$0 \fun \cI_W^\vee[1] \fun E \fun \cI_Z(H) \fun 0$$
where $Z,W \subset S$ are zero-dimensional subschemes of the same length.

\medskip

Moreover, if $E$ is one of the objects above, and $E$ is not $t$-stable for some $t > \frac 16$, then
$E$ is destabilized by a twisted ideal sheaf $\cI_Y(H) \subset E$ for some zero-dimensional subscheme $Y \subset S$.
\end{thm}

\pf By Proposition \ref{Stable}, any $t$-stable object is either of the form $i_*L_C$ or else fits in a sequence:
$$0 \fun K \fun E \fun Q \fun 0$$
where $Q$ is an $H$-stable torsion-free sheaf of odd rank $2n+1$ and $c_1 = (n+1)H$ and $K$ is the shifted derived dual of a torsion-free sheaf. But by the Bogomolov-Gieseker inequality, 
$$\ch_2(Q) \le \frac{(n+1)^2H^2}{2(2n+1)}$$
and then
$$\mbox{Re}(Z_t(Q)) \ge \frac {H^2}2\left(t^2(2n+1)  - \frac 1{4(2n+1)} \right)$$
Thus if $t \ge \frac 12$, there are no such $t$-stable objects (so the only $t$-stable objects are
of the form $i_*L_C$), if $t \ge \frac 16$, then there are none such with $r = (2n+1) \ge 3$, etc.

\medskip

The last part of the theorem now follows from 
Lemmas \ref{StableL} and \ref{StableExt}. 

\qed

\section{Families and Walls}

For $t > \frac 12$, the moduli of $t$-stable objects of $\cA^\#$ is the moduli space:
$$\cM := \cM_S\left(0,H,\frac {H^2}2\right)$$
of (Gieseker-stable) {\it coherent sheaves} on $S$ of the form $i_*L_C$. As $t$ crosses the critical values:
$$\frac 12, \sqrt{\frac 14 - \frac 2{H^2}}, \sqrt{\frac 14 - 2\frac {2}{H^2}}, \cdots > \frac 16$$
the $t$-stability changes, as subobjects of certain coherent sheaves $i_*L_C$  (or more generally,
objects of $\cA^\#$ from Lemma \ref{StableExt}) of the form $\cI_{Z}(H)$
achieve zero (and then positive) slope. 
The resulting birational modifications of $\cM$ as $t$ passes over critical points 
can be predicted, but are only carried out (in \S 5)
in case $S$ is $K$-trivial, because it is only in that case that we can prove that the desired birational
transformations (which are then Mukai flops) actually exist.

\medskip

\begin{defi} For (quasi-projective) schemes $X$, the objects $\cE_X$ of the bounded derived category
${\cD}(S\times X)$ are {\it families of objects} of $\cD(S)$ parametrized by $X$.
\end{defi}

\begin{defi} A family $\cE_X$ is a {\it flat family of objects of $\cA^\#$} if
the (derived) restrictions to the fibers:
$$\cE_x := Li_{S\times x}^*{\cE}$$
are objects of $\cA^\#$ for all closed points $x\in X$ (via the isomorphism $S \times x \cong S$).
\end{defi}

\begin{rem} This is a good analogue of the flat families of coherent sheaves on $S$.
The category of flat families of objects of $\cA^\#$, like the category of flat families
of coherent sheaves, is not abelian, but Abramovich-Polishchuk define an analogue of the abelian category of coherent sheaves on $X\times S$ \cite{AP06} (at least in the case where $X$ is smooth).
We will not need to make use of this abelian category.
\end{rem}

\begin{ex} The universal family of coherent sheaves:
$$\cU \rightarrow S \times \cM$$
for the moduli space $\cM = \cM_S\left(0,H,H^2/2\right)$ is a flat family of objects of $\cA^\#$ (all coherent sheaves!).
\end{ex}

\begin{ex} Let $S[d]$ be the Hilbert scheme of length $d$ subschemes of $S$, with universal subscheme:
$$\cZ \subset S \times S[d]$$
Then the sheaf $\cI_\cZ(H) := \cI_\cZ \otimes \pi_1^*\cO_S(H)$ is a flat family of 
objects of $\cA^\#$ (and of coherent sheaves).
\end{ex}

\begin{ex} The shifted derived dual $\cI_\cZ^\vee[1]$ (in $\cD(S\times S[d])$)  is a flat family
of objects of $\cA^\#$ . 

\medskip

Indeed, it is a consequence of the flatness of the coherent sheaf $\cI_\cZ$ over $S[d]$ that:
$$Li_{S\times \{Z\}}^*\cI_\cZ^\vee[1]  = \cI_Z^\vee[1]$$
for each $Z\in S[d]$.

\end{ex}

Our goal is to produce flat families of objects of $\cA^\#$ parametrizing extensions of the form:
$$0 \fun \cI_Z(H) \fun E \fun \cI_W^\vee[1] \fun 0 \ \ \mbox{and}\ \ 0 \fun \cI_W^\vee[1] \fun E \fun \cI_Z(H) \fun 0$$
that are exchanged under the wall-crossing. These will both be {\bf projective bundles} 
when $S$ is $K$-trivial, thanks to the 
following vanishing result:

\begin{pro}\label{Van} Let $S$ be a smooth surface with Pic$(S) = \bZ[H]$. Then:
$$\mbox{H}^i(S,\cI_Z\otimes \cI_W\otimes \cO_S(H + K_S)) = 0; \ \ \mbox{for}\ \ i = 1,2$$
for all subschemes $Z,W \subset S$ of (the same) length $d$ provided that:
$$d < \frac {H^2}8$$
\end{pro}

\begin{rem} In the case $d = 1$, this is a weak form of Reider's Theorem \cite{Fri98}, 
since it amounts to saying that
$H + K_S$ is very ample if $H$ is ample, generating Pic$(S)$, and $H^2 \ge 9$.
\end{rem}

\pf In the derived category $\cD(S)$,
$$H^i(S,\cI_Z\otimes \cI_W\otimes \cO_S(H + K_S)) \cong \Ext^{i+1}_{\cD(S)}(\cI_W^\vee[1],\cI_Z(H+K_S))$$
and by Grothendieck duality:
$$\Ext^{i+1}_{\cD(S)}(\cI_W^\vee[1],\cI_Z(H+K_S)) \cong \Ext^{1-i}_{\cA^\#}(\cI_Z(H),\cI_W^\vee[1])^*$$
This immediately gives $H^2(S,\cI_Z\otimes \cI_W\otimes \cO_S(H + K_S)) = 0$ (which was easy to check anyway),
and identifies the space of extensions of $\cI_Z(H)$ by $\cI_W^\vee[1]$  with  
$H^0(S,\cI_Z\otimes \cI_W\otimes \cO_S(H + K_S))^*$. But it also identifies:
$$H^1(S,\cI_Z\otimes \cI_W\otimes \cO_S(H + K_S))^* \cong \Hom_{\cA^\#}(\cI_Z(H),\cI_W^\vee[1])$$
and we may conclude that this is zero if we can find a value of $t > 0$ such that:
$$\cI_Z(H), \cI_W^\vee[1] \ \mbox{are both $t$-stable and}\ \ \mu_t(\cI_Z(H)) > \mu_t(\cI_W^\vee[1])$$
But the ``wall'' where these two slopes coincide is precisely at:
$$t = \sqrt{\frac 14 - d\frac 2{H^2}}$$
which satisfies $t > 0$ when $d < \frac {H^2}8$, as desired.

\medskip

\begin{rem} It is  interesting that this vanishing theorem can be proved using stability conditions, 
following more or less immediately from the Hodge Index Theorem and the Bogomolov-Gieseker inequality. 
A stronger inequality valid for K3 surfaces will give a stronger vanishing result in 
\S 6.

\end{rem}

\begin{pro}\label{Pd}  Assume the vanishing of Proposition \ref{Van}. If $K_S \ge 0$, the projective bundle:
$$\bP_d \rightarrow S[d]\times S[d] \ \mbox{with fibers}\ \ \bP(H^0(S,\cI_Z\otimes \cI_W(H+K_S)))$$
supports a universal family $\cE_d$
(on $S\times \bP_d$) of extensions of objects of $\cA^\#$ of the form:
$$0 \fun \cI_Z(H+ K_S) \fun E \fun \cI_W^\vee[1] \fun 0$$
For any $S$, the dual projective bundle $\bP_d^\vee$ supports a universal family $\cF_d$ of extensions 
of the form:
$$0 \fun \cI_W^\vee[1] \fun F \fun \cI_Z(H) \fun 0$$
\end{pro}

\pf Let:
$$\cZ_{12}, \cZ_{13} \subset S \times S[d] \times S[d]$$
be the pull-backs of $\cZ \subset S\times S[d]$ via the projections:
$\pi_{12}, \pi_{13}: S \times S[d] \times S[d] \rightarrow S \times S[d]$
and consider the (a priori derived) object:
$$RHom(\cI_{\cZ_{13}}^\vee[1], \cI_{\cZ_{12}}(H+K_S))[1] \cong \cI_{\cZ_{13}} \stackrel
L \otimes \cI_{\cZ_{12}}(H+K_S)$$

Since the ideal sheaf $\cI_\cZ$ admits a two-step resolution by vector bundles (see \cite{ES98})
it follows that 
$\cI_{\cZ_{13}} \stackrel L \otimes \cI_{\cZ_{12}}(H)$ is 
(equivalent to) a flat coherent sheaf over $S[d]\times S[d]$.

\medskip
 
 Let $\pi: S\times S[d]\times S[d] \rightarrow S[d]\times S[d]$ be the projection. We may set:
 $$\bP_d := \bP\left(\pi_*\left(\cI_{\cZ_{13}}\otimes \cI_{\cZ_{12}}(H+K_S)\right)\right)$$
 since the push-forward is locally free, by base change.
  
\medskip

Next, we turn to the construction of the universal family on $S\times \bP_d$. This 
should morally be thought of as a single extension of the form:
$$0 \fun \rho^*\cI_{\cZ_{12}}(H + K_S) \fun \cE_d \fun \rho^*\cI_{\cZ_{13}}^\vee[1] 
\otimes \cO_{\bP_d}(-1)\fun 0$$
where
$$\rho: S\times \bP_d \rightarrow S\times S[d] \times S[d]$$
is the projection. But we have avoided any mention of an abelian category of
objects on $S\times \bP_d$ containing both
$\rho^*\cI_{\cZ_{12}}(H+K_S)$ and 
$\rho^*\cI_{\cZ_{13}}^\vee[1] \otimes \cO_{\bP_d}(-1)$. Instead, we make the 
construction using distinguished triangles. There is a canonical element:
$$\id \in \Gamma(S\times \bP_d, \rho^*\left(\cI_{\cZ_{13}}\otimes \cI_{\cZ_{12}}(H+K_S)\right) \otimes \cO_{\bP_d}(1))$$
$$ = \Gamma(S\times S[d]\times S[d], \cI_{\cZ_{13}}\otimes \cI_{\cZ_{12}}(H+K_S) \otimes \rho_* \cO_{\bP_d}(1))$$
$$= \Gamma(S[d]\times S[d], \pi_*(\cI_{\cZ_{13}}\otimes \cI_{\cZ_{12}}(H+K_S)) \otimes \pi_*
(\cI_{\cZ_{13}}\otimes \cI_{\cZ_{12}}(H+K_S))^*)$$
which can be alternatively thought of as the canonical element:
$$f_\id \in RHom_{S\times \bP_d}(\rho^*\cI_{\cZ_{13}}^\vee \otimes \cO_{\bP_d}(-1), \rho^*
\cI_{\cZ_{12}}(H+K_S) )$$

With this canonical element, we may form the cone and distinguished triangle:
$$(**) \ \ \ \cdots \fun \rho^*\cI_{\cZ_{13}}^\vee \otimes \cO_{\bP_d}(-1)
\stackrel {f_\id} \fun \rho^* \cI_{\cZ_{12}}(H+K_S)  
\fun \cE_d \fun \rho^*\cI_{\cZ_{13}}^\vee[1]  \otimes \cO_{\bP_d}(-1)
\stackrel {f_\id[1]} \fun \cdots $$

If $K_S \ge 0$, then this ``universal'' distinguished triangle has the property that each:
$$Li_{S\times \epsilon}^*(**): \cdots \fun \cI_Z(H+K_S) \fun \cE_d|_{S\times \epsilon} \fun 
\cI_W^\vee[1] \fun \cdots$$
{\bf is} the short exact sequence (in $\cA^\#$) corresponding to the extension (modulo scalars):
$$\epsilon \in \bP(H^0(S,\cI_W\otimes \cI_Z(H + K_S))) \cong \bP(\Ext^1_{\cA^\#}(\cI_W^\vee[1],\cI_Z(H+K_S)))$$

Turning next to the family $\cF_d$, we define, similarly,
$$\rho^\vee: S \times \bP_d^\vee \fun S\times S[d]\times S[d]$$
and as above, the key point is the existence of a canonical morphism:
$$f^\vee_{id}:   (\rho^\vee)^*\cI_{\cZ_{12}}(H)[-1] \fun (\rho^\vee)^*\cI_{\cZ_{13}}^\vee[1] \otimes \cO_{\bP_d^\vee}(1)$$
which will, in turn, define the distinguished triangle:
$$\cdots \fun (\rho^\vee)^*\cI_{\cZ_{12}}(H)[-1] \fun (\rho^\vee)^*\cI_{\cZ_{13}}^\vee[1] \otimes \cO_{\bP_d^\vee}(1) \fun \cF_d \fun (\rho^\vee)^*\cI_{\cZ_{12}}(H) \fun \cdots$$
which is the desired universal family (whether or not $K_S \ge 0$!)

\medskip

But this canonical morphism is obtained from Serre duality \cite{Cal05}:
$$RHom((\rho^\vee)^*\cI_{\cZ_{12}}(H)[-1], (\rho^\vee)^*\cI_{\cZ_{13}}^\vee[1] \otimes \cO_{\bP_d^\vee}(1))$$
$$ \cong RHom((\rho^\vee)^*\left(\cI_{\cZ_{12}} 
\otimes \cI_{\cZ_{13}} (H)\right), \cO_{\bP_d^\vee}(1)[2])$$
$$\cong RHom((\rho^\vee)^*\left(\cI_{\cZ_{12}} 
\otimes \cI_{\cZ_{13}} (H+K_S)\right), \cO_{\bP_d^\vee}(1)\otimes \cO_S(K_S)[2])$$
$$ = 
RHom((\rho^\vee)^*\left(\cI_{\cZ_{12}} 
\otimes \cI_{\cZ_{13}} (H+K_S)\right), \pi^! \cO_{\bP_d^\vee}(1))$$
$$ \cong 
RHom\left(\pi_*\left(\cI_{\cZ_{12}} 
\otimes \cI_{\cZ_{13}} (H+K_S)\right), \pi_*\left(\cI_{\cZ_{12}} 
\otimes \cI_{\cZ_{13}} (H+K_S)\right)\right)$$
\qed

\section{Mukai Flops} 

In this section, and for the rest of the paper, we assume that $K_S = 0$, for the following reason.

\medskip

Moduli spaces $\cM = \cM_S(r,c_1,\emph{ch}_2)$ of $H$-stable coherent sheaves on a $K$-trivial surface are symplectic, 
meaning that there is a skew-symmetric 
isomorphism on the tangent bundle:
$$\omega : T\cM \fun T^*\cM$$
The form is given by the natural isomorphism of Serre duality (see \cite{Muk84}): 
$$\Ext^1_{\cO_S}(E,E) \cong \Ext^1_{\cO_S}(E,E)^*$$
Note that this argument could apply as well to moduli spaces of stable objects of $\cA^\#$, taking
$\Ext^1_{\cA^\#}(E,E)$, once moduli spaces are shown to exist(!)

\medskip

Varieties with a symplectic structure 
are necessarily even-dimensional. When such a variety is 
equipped with an appropriate ``Lagrangian'' subvariety, then it always admits an
elementary birational transformation (nowadays known as a {\it Mukai flop}):

\medskip

\begin{thm2} [Theorem 0.7 of  \cite{Muk84}] Let $X$ be a symplectic variety, and 
let $P$ be a $\bP^n$-bundle, over a base $B$, contained in $X$ in codimension $n \ge 2$. Then there is a  birational map, denoted elm$_P:X--\!\!\!>X'$ with the following properties:

\medskip

1) $X'$ contains the dual $\bP^n$ bundle $P'$ over $B$ and has a symplectic structure $\omega'$
which coincides with $\omega$ outside of $P'$

\medskip

2) elm$_P$ is the composite of the blowing up $\sigma^{-1}:X --> \widetilde X$ along $P$ and the 
blowing down $\sigma': \widetilde X \fun X'$ of the exceptional divisor $D = \sigma^{-1}(P)$ onto $P'$. Moreover, $D \subset P \times_B P'$ is the two-step flag bundle over $B$, and $\cO_{\widetilde X}(D)|_D \cong \cO_{P\times_B P'}(-1,-1)|_D$.

\end{thm2}

\medskip

\begin{thm}\label{moduli} Fix $S$ with $K_S = 0$ and Pic$(S) = \bZ[H]$, (i.e. $S$ is a K3 surface) and let 
$$\left\{t_d = \sqrt{\frac 14 - d\frac 2{H^2}} \ | \ d = 0,1,2,... < \frac{H^2}8\right\}$$
be the set of ``rank one'' critical values for stability conditions $(Z_t,\cA^\#)$. Then
for each $t > \frac 16$ and away from the critical set, there is a smooth, proper moduli space:
$$\cM_t := \cM_t\left(0,H,\frac {H^2}2\right)$$
which together with a suitable {\bf coherent sheaf} $\  \cU_t$ on $S\times \cM_t$ represents the functor:
 isomorphism classes of flat families of $t$-stable objects.

\end{thm}

{\bf Proof:} The moduli space $\cM_t$ for $t > \frac 12$ is 
the ``classical'' space $\cM_S(0,H,\frac{H^2}2)$ of rank one torsion-free sheaves of degree $H^2$ 
on curves $C\in |H|$. This admits a universal coherent sheaf by 
Geometric Invariant Theory. The general
proof consists of three parts, which carry out an induction that constructs 
each $\cM_{t_d - \epsilon}$ out of $\cM_{t_d + \epsilon} (= \cM_{t_{d-1} - \epsilon})$ 
near each of the ``walls'' $t_d  > \frac 16$. 

\medskip

Assume that $t_d > \frac 16$ and $\cM_{t_d  + \epsilon}$, together with ``universal'' coherent sheaf $\cU_{t_d +\epsilon}$ 
on $S\times \cM_{t_d + \epsilon}$, is the smooth, proper moduli space representing the functor 
of isomorphism classes of $t_d + \epsilon$-stable objects with the given invariants. Then:

\medskip

\nt {\bf Step 1:} There is a natural embedding ${\bP}_d \subset  \cM_{t_d + \epsilon}$ of 
the projective bundle from Theorem \ref{Pd}  that parametrizes all the objects of 
$\cM_{t_d + \epsilon}$ 
that are not $t_d -\epsilon$-stable. 

\medskip

\nt {\bf Interlude:} Construct  ${\bP}^\vee_d \subset \cM'$  as the Mukai flop of ${\bP}_d \subset 
\cM := \cM_{t_d + \epsilon}$.

\medskip

\nt {\bf Step 2:} There  is a coherent sheaf $\cU'$ on $S\times \cM'$ naturally obtained as the 
``Radon transform'' across the Mukai flop of the universal coherent sheaf $\cU_{t_d + \epsilon}$ on 
$S\times \cM_{t_d + \epsilon}$, such that:

\medskip

\nt {\bf Step 3:} $\cM'$ together with the sheaf $\cU'$ is the desired 
$\cM_{t_d - \epsilon}$ (and family $\cU_{t_d - \epsilon}$).

\medskip 

{\bf Proof of Step 1:} If $Z,W \subset S$ are subschemes of length $d$, consider an object 
$E$ of $\cA^\#$ given as an extension:
$$0 \rightarrow \cI_Z(H) \rightarrow E \rightarrow \cI_W^\vee[1] \rightarrow 0$$
First, recall that both
$\cI_Z(H)$ and $\cI_W^\vee[1]$ are $t$-stable (Lemma \ref{LambdaI}). Since:
$$\Real(Z_t(\cI_Z(H))) =  d + \frac{H^2}2(t^2 - \frac 14) \ \mbox{and}\ 
\Real(Z_t(\cI_W^\vee[1])) = - \Real(Z_t(\cI_Z(H)))$$
it  follows that $E$ is not $t$-stable if $t \le t_d$ (recall that $t_d$ solves $d + \frac{H^2}2(t_d^2 - \frac 14) = 0$).
But we claim that $E$ is $t$-stable for $t_d < t < t_{d-1}$. To see 
this, consider the cohomology sequence of sheaves associated to the extension defining $E$:
$$0 \fun H^{-1}(E) \fun \cO_S \fun \cI_Z(H) \fun H^0(E) \fun T \fun 0$$
From this it follows that either $H^{-1}(E) = 0$, in which case $E = H^0(E) =  i_*L_C$ for some  (necessarily torsion-free!) rank one sheaf on $C$, 
or else $H^{-1}(E) = \cO_S$, and then $E$ fits in an extension: 
$$0 \fun \cI_{W'}^\vee[1] \fun E \fun \cI_{Z'}(H) \fun 0$$
where $\cI_{Z'}(H) = H^0(E)/tors(H^0(E))$. Moreover, this sequence only splits if $Z = Z'$ and $W = W'$, and the original exact sequence is split. Thus the second part of 
Theorem \ref{tstable} applies, and if $E$ were not $t$ stable, 
then it would be destabilized by 
an ideal sheaf $\cI_Y(H) \subset E$. On the other hand, such ideal sheaves satisfy $\Real(Z_t(\cI_Y(H))) = d' + \frac{H^2}2(t^2 - \frac 14)$, where $d' = len(Y)$, so if $E$ were destabilized by such a sheaf, and if 
$t_d < t < t_{d-1}$, then $len(Y) \le d-1$, the induced map  
$\cI_Y(H) \fun \cI_W^\vee[1]$ is the zero map (otherwise it would destabilize $\cI_W^\vee[1]$!), and so $\cI_Y(H) \subset \cI_Z(H)$, which contradicts the fact that $len(Y) < len(Z)$.

\medskip

Thus the non-split extensions parametrized by ${\bP}_d$ produce $t$-stable objects (for $t_d < t < t_{d-1}$) and
the family of Theorem \ref{Pd} defines a morphism:
$$i_d: {\bP}_d \rightarrow \cM_{t}$$
We claim that $i_d$ is an embedding. First, if $i_d(\epsilon) = i_d(\epsilon')$ , where $\epsilon, \epsilon'$
are extension classes defining isomorphic objects $E,E'$, 
then because $E$ and $E'$ are both $t$-stable, 
it follows that the isomorphism is a multiple of the 
identity map. Moreover, since Hom$(\cI_Z(H),\cI_{W'}^\vee[1]) = 0$ for all 
$Z,W$ satisfying $len(Z) = len(W') = d$, it follows that the isomorphism $E \cong E'$ induces 
vertical isomorphisms in the following diagram:
$$\begin{array}{cccccccccccccccc}
\epsilon: & 0 & \fun & \cI_Z(H) & \fun & E & \fun & \cI_W^\vee[1] & \fun & 0\\
& & & \|\wr & & \|\wr & & \|\wr \\
\epsilon': & 0 & \fun & \cI_{Z'}(H) & \fun & E'  & \fun & \cI_{W'}^\vee[1]  & \fun & 0\end{array}$$
{\it all of which are multiples of the identity}. Thus
$Z = Z', W = W'$ and $\epsilon = \lambda \epsilon'$ for some non-zero scalar $\lambda$.
In other words, the equivalence classes of the extensions modulo scalars satisfy
$[\epsilon] = [\epsilon']$ (in $\bP_d$). Thus $i_d$ is injective.

\medskip

To complete the proof that $i_d$ is an embedding, we need to study the induced map on tangent spaces. The tangent space to $\cM_t$ at a point $E\in \cM_t$ is easiest to describe. It is:
$$\Ext^1_{\cA^\#}(E,E)$$
(the same as the tangent space to the stack...see Step 3). If $E = i_d([\epsilon])$, where:
$$\epsilon: 0 \fun \cI_Z(H) \fun E \fun \cI_W^\vee[1] \fun 0$$
is a (non-split) extension, then $\Ext^1_{\cA^\#}(E,E)$ fits into a long exact sequence:
$$0 \fun V  \fun \Ext^1_{\cA^\#}(E,E) \fun  \Ext^1_{\cA^\#}(\cI_Z(H), \cI_W^\vee[1]) \stackrel {\epsilon^\vee} \fun 
\Ext^2_{\cA^\#}(\cI_W^\vee[1], \cI_W^\vee[1]) = \bC \fun 0$$
where $V$ is identified with the tangent space to $\bP_d$ at the point $[\epsilon]$ via:
$$ \Ext^1_{\cA^\#}(\cI_Z(H),\cI_Z(H)) \cong T_{S[d]}(Z), \ \ 
 \Ext^1_{\cA^\#}(\cI_W^\vee[1],\cI_W^\vee[1]) \cong T_{S[d]}(W)$$
 and 
$$0 \fun \bC = \Hom(\cI_Z(H),\cI_Z(H)) \stackrel \epsilon \fun \Ext^1_{\cA^\#}(\cI_W^\vee[1],\cI_Z(H)) \fun
V \fun$$
$$ \fun \Ext^1_{\cA^\#}(\cI_Z(H),\cI_Z(H)) \oplus \Ext^1_{\cA^\#}(\cI_W^\vee[1],\cI_W^\vee[1])  \fun 0$$

A standard deformation theory argument shows that the induced map 
from $V \cong T_{\bP_d}([\epsilon])$ to $T_{\cM}(E)$ is the differential of $i_d$. 
Thus $i_d$ is an embedding, and the normal space at $E$ is naturally identified 
with the kernel of the map:
$$\Ext^1_{\cA^\#}(\cI_Z(H),\cI_W^\vee[1]) \cong \Ext^1(\cI_W^\vee[1],\cI_Z(H))^* \stackrel{\epsilon ^\vee}
\fun \bC$$

\medskip

\nt {\bf Interlude:} By induction (or else directly!), $\cM_t$ is a symplectic variety,  and:

\medskip

$\bullet$ \ $\dim(\cM_t) = \dim(\Ext^1_{\cD(S)}(E,E)) = 2 + H^2$

\medskip

$\bullet$ \ $\dim(\bP_d) = 2\dim(S[d]) + (\chi(S,\cI_Z\otimes \cI_W\otimes \cO_S(H)) - 1) = 4d + (1 + \frac{H^2}2 - 2d)$

\medskip

\nt so that, indeed, the embedding of Step One satisfies:
$$\codim(\bP_d) = \frac{H^2}2 + 1 - 2d = \mbox{fiber dimension}$$
and there is a Mukai flop: $ \bP_d  = P \subset \cM := \cM_{t_d + \epsilon} --> \cM' \supset P'$. 
We can now describe all the points of $\cM'$:

\medskip

$\bullet$ The points of $\cM' - P' = \cM - P$ correspond to all the objects of $\cA^\#$ with the given 
invariants that are both $t_d + \epsilon$-stable and $t_d - \epsilon$-stable.

\medskip

$\bullet$ Via the isomorphism:
$\Ext^1_{\cA^\#}(\cI_Z(H), \cI_W^\vee[1]) \cong \Ext^1_{\cA^\#}(\cI_W^\vee[1],\cI_Z(H))^*$
the points of $P' = \bP_d^\vee$ correspond to non-zero extensions (modulo scalars) of the form:
$$(*) \ \ 0 \fun \cI_W^\vee[1] \fun E \fun \cI_Z(H) \fun 0$$
with $len(W) = len(Z) = d$.

\medskip

It is evident that such extensions define objects 
$E$ of $\cA^\#$ that are not $t_d + \epsilon$-stable. 
On the other hand, if $E$ is an object of $\cA^\#$ with the given invariants 
that is $t_d - \epsilon$ stable, then by Theorem \ref{tstable}, $E$ is either of the form
$i_*L_C$ or else is an extension of the form $(*)$ with $len(W) = len(Z) \le d$.
By the second part of that theorem, any sheaf $i_*L_C$ 
or extension of the form $(*)$ with 
$len(W) = len(Z) < d$ that is $t_d - \epsilon$ stable is also $t_d + \epsilon$ stable.
This just leaves the points of $P'$, which are all $t_d - \epsilon$-stable, by the same 
argument as in Step 1. Thus the points of $\cM'$ are in bijection with the 
isomorphism classes of all the $t_d - \epsilon$ stable objects of $\cA^\#$ (with the 
given invariants)!
 
\medskip

\nt {\bf Proof of Step 2:} Let $\sigma: \widetilde \cM \fun \cM$ be the blow-up along $\bP_d$, and let  
$$\widetilde \cU := L\sigma^*\cU$$
be the (derived) pullback of the coherent sheaf $\cU$ from $S\times \cM$ to $S\times \widetilde \cM$.
We will prove that $\widetilde \cU$ is in fact a coherent sheaf, and then construct $\cU'$ by 
descending a (generalized) elementary modification of $\widetilde \cU$. Carrying out this 
elementary modification will require two(!) applications of the octahedral axiom.
 
 \medskip
 
\begin{lem}\label{Sheaf} 
Suppose $\cE_X$ is a flat family of objects of $\cA^\#$ over a connected 
quasi-projective base $X$.
If $\cE_x$ is (quasi-isomorphic to) a coherent sheaf on $S$ for some closed point
$x\in X$, then $\cE_X$
is (quasi-isomorphic to) a coherent sheaf on $S\times X$.
\end{lem}

\pf The object $\cE_X$ can be represented (in $\cD^b(S\times X)$) by a two-term complex
$[A \fun B]$. Moreover, by pulling back under a surjective map $V \fun B$, we may assume,
without loss of generality, that $B = V$ is a locally free sheaf. Since $\cE_X$ is a flat family
of objects of $\cA^\#$, in particular we know that each restriction $\cE_x$ has cohomology 
only in degrees $-1$ and $0$. Thus it follows that $A$ is {\it flat} as a coherent sheaf over $X$.
Now suppose that some $\cE_x$ is a coherent sheaf on $S$. Then the kernel of the map 
$A \fun V$ must be zero generically over $X$, and, if non-zero, would determine an embedded point
of $A$ that does not dominate $X$. Such coherent sheaves are not flat over $X$. Thus
$A \fun V$ is injective, and $\cE_X$ is (quasi-isomorphic to) a coherent sheaf.

\qed

\begin{cor}\label{PdFamily}
The flat families $\cE_d$  and $Li^*\cU$ on $S\times \bP_d$ are (quasi-isomorphic to) coherent sheaves.
\end{cor}

\pf Among the extensions parametrized by $\cE_d$ are the extensions yielding:
$$0 \fun \cO_S \fun \cI_Z(H) \fun \cO_C(H + W - Z) \fun \cO_W \fun 0$$
where $C\in |H|$ is a smooth curve, and $W,Z \subset C$ are disjoint effective divisors of degree $d$ 
on $C$
(see the example preceding Lemma \ref{StableExt}).
Such an extension $\epsilon$ defines $\cO_C(H+W-Z)$ as $(\cE_d)_{\epsilon}$.
\qed

\medskip

Consider now the pair of coherent sheaves on $S\times \bP_d$:
$$\cU|_{S\times \bP_d} \ \mbox{and}\ \cE_d$$
We cannot conclude that the two sheaves are isomorphic, because there is a built-in 
ambiguity from the fact that $\cU$ and $\cU\otimes \cL_{\cM}$ 
give equivalent universal families for any line bundle $\cL_{\cM}$ on $\cM$ (this is, as in the case of the Jacobian, the only ambiguity).
But we can ``match'' them as closely as we need with the following:

\medskip

\begin{lem}\label{Match}
The map $\emph{Pic}(\cM) \fun \bZ = \emph{Pic}(\bP_d)/\emph{Pic}(S[d]\times S[d])$ is surjective.
\end{lem}

\pf For fixed disjoint reduced subschemes $W,Z \subset S$ of length $d$, the fiber of $\bP_d$ over 
the points $(Z,W)$ is, outside a subvariety of codimension $> 1$, isomorphic to the 
linear series $\bP^{g-2d} = |\cO_S(H) \otimes \cI_Z \otimes \cI_W|$ of curves passing through the 
points of $Z$ and $W$ (determining the line bundle $\cO_C(H+W-Z)$ as in the proof of 
the Corollary above). The line bundle $\pi^*\cO_{\bP^g}(1)$ pulled back from 
the ``linear series map'' $\pi:\cM_S(0,H,g-1) \fun \bP^g$ carries over to a line bundle
on $\cM$ (across all previous Mukai flops), which agrees with $\cO_{\bP^{g-2d}}(1)$ off
codimension $>1$. Thus this line bundle on $\cM$ generates the relative Picard group
of $\bP_d$ over $S[d]\times S[d]$, as desired.
\qed

\begin{cor}\label{Restrict}
There is a choice of ``Poincar\'e'' coherent sheaf $\cU$ on $S\times \cM$ and line bundle 
$\cL$ on $S[d]\times S[d]$ such that $\cU|_{S\times \bP_d} \cong \cE_d\otimes \cL$, hence $\cU|_{S\times \bP_d}$ fits in a
distinguished triangle of the following form:
$$\cdots \fun \rho^* \left(\cI_{\cZ_{12}}(H)\otimes \cL\right) 
\fun \cU|_{S\times \bP_d} \stackrel u \fun \rho^*\left(\cI_{\cZ_{13}}^\vee[1]  \otimes \cL\right)
\otimes \cO_{\bP_d}(-1)  \fun \cdots $$
($\rho: S\times \bP_d \fun S\times S[d]\times S[d]$ and other notation from the 
proof of Proposition \ref{Pd}).
\end{cor}

\pf By the Lemma, we can match $\cU|_{S\times \bP_d}$ and $\cE_d$ up to the twist of a line bundle 
$\cL$ pulled
back from $S[d]\times S[d]$. The distinguished triangle is then the corresponding twist
of the distinguished triangle defining $\cE_d$ in the proof of Proposition \ref{Pd}.
\qed

\medskip

Now, let $i_D:D \hookrightarrow \widetilde \cM$ be the exceptional divisor of the blow-up, and let 
$D_S = S\times D \subset S\times \widetilde \cM$ with projections $p :D_S \fun 
S\times \bP_d$ and $p':D_S \fun S\times \bP_d^\vee$.
Then we define a 
$\widetilde \cU' \in \cD(S\times \widetilde \cM)$ via the distinguished triangle:
$$\cdots \fun \widetilde \cU' \fun \widetilde \cU \stackrel {u\circ r} \fun {i_{D_S}}_*p^*\left(\rho^*\left(\cI_{\cZ_{13}}^\vee[1]  \otimes \cL\right) \otimes \cO_{\bP_d}(-1)\right) \fun \widetilde \cU'[1] \fun \cdots$$
where the ``restriction map'' $r$ fits in the distinguished triangle (of coherent sheaves!):
$$\cdots \fun \widetilde \cU(-D_S) \fun \widetilde \cU 
\stackrel r \fun {i_{D_S}}_*\widetilde \cU|_{D_S} \fun \widetilde \cU(-D_S)[1] \fun\cdots $$
and $u$ is (by abuse of notation) the map from the distinguished triangle of Corollary \ref{Restrict},
pulled back to $D_S$ and pushed forward to $S\times \widetilde \cM$. The octahedral 
axiom applied to the morphisms $u$ and $r$ now produces 
a distinguished triangle:
$$\cdots \fun \widetilde \cU(-D_S) \fun \widetilde \cU' 
 \stackrel v \fun {i_{D_S}}_*p^*\rho^* \left(\cI_{\cZ_{12}}(H)\otimes \cL\right) \fun \widetilde \cU(-D_S)[1] \fun \cdots $$

\rem This is the derived category version of an elementary modification of a coherent sheaf by a quotient
sheaf supported on a divisor. Although the object we are modifying, $\widetilde \cU$, is indeed a coherent sheaf, it is being modified by a ``quotient'' which is not a coherent sheaf. Nevertheless, we 
constructed the modified object $\widetilde \cU'$ 
as the (shifted) cone of the morphism to the ``quotient'' object supported on $D_S$.

\medskip

Next, we claim that the (derived) restriction of $\widetilde \cU'$ to $D_S$ satisfies:
$$ L{i^*_{D_S}} \widetilde \cU' \cong p'^*\cF_d\otimes \cL$$
where $\cF_d$ is the ``universal extension'' on $S\times \bP_d^\vee$ from Proposition
\ref{Pd}:
$$\cdots  \fun (\rho^\vee)^*\left(\cI_{\cZ_{13}}^\vee[1]  \right)
\otimes \cO_{\bP_d^\vee}(+1)  \fun \cF_d  \fun  (\rho^\vee)^* \left(\cI_{\cZ_{12}}(H) \right) 
\fun \cdots$$

Note that $\cF_d$ (and its pullback) is definitely {\bf not} a coherent sheaf, although 
$\widetilde \cU'$ itself {\it will} be.
We see the claim with another application of the octahedral axiom, to the two morphisms:
$\widetilde \cU' \stackrel r \fun {i_{D_S}}_* L{i^*_{D_S}}\widetilde \cU' $
and the push-forward of
$$L{i^*_{D_S}}\widetilde \cU' \stackrel v \fun p^*\rho^*\left(\cI_{\cZ_{12}}(H)\otimes \cL\right)
\cong  p'^*(\rho^\vee)^*\left(\cI_{\cZ_{12}}(H)\otimes \cL\right)$$
($v$ is the left adjoint of the map $v$ defined by the first application of the octahedral axiom!)
Let $K$ be defined by the distinguished triangle:
$$ \cdots \fun K \fun L{i^*_{D_S}}\widetilde \cU' \stackrel v \fun p^*\rho^*\left(\cI_{\cZ_{12}}(H)\otimes \cL\right) \fun K[1] \fun \cdots $$
Then the octahedral axiom applied to $r$ and ${i_{D_S}}_*v$ gives:
$$\cdots \fun \widetilde \cU'(-D_S) \fun \widetilde \cU(-D_S) \fun {i_{D_S}}_*K \fun \cdots$$
which in turn allows us to conclude that, as desired:
$$K \cong p^*\left(\rho^*\left(\cI_{\cZ_{13}}^\vee[1] \otimes \cL\right)
\otimes \cO_{\bP_d}(-1)  \right) \otimes \cO_{D_S}(-D_S)$$
$$ \cong 
p'^*\left((\rho^\vee)^*\left(\cI_{\cZ_{13}}^\vee[1] \otimes \cL\right)
\otimes \cO_{\bP_d^\vee}(+1)  \right)$$

\medskip

\nt {\bf Main Point of Step Two:} There is a coherent sheaf  $\cU'$ on $S\times \cM'$ 
which is a flat family of objects of $\cA^\#$ such that
the restrictions $\cU'|_{S\times \{m'\}}$ are in bijection with the set of all $t_d-\epsilon$-stable objects 
of $\cA^\#$ (with the given invariants). 
This coherent sheaf is obtained by descending the object $\widetilde \cU'\in \cD(S\times \widetilde \cM)$, defined above, to $\cD(S\times \cM')$. That is:
$$L_{\sigma'}^*\cU' \cong \widetilde \cU'$$

\pf 
Since $\widetilde \cU' $ and $\widetilde \cU$ coincide away from the exceptional divisor $D_S$, and 
$L{i_{D_S}}^*\widetilde \cU'$ is the pullback of a 
(universal) family of 
extensions of objects of $\cA^\#$, it follows that $\widetilde \cU'$ is a flat family of objects of $\cA^\#$. 
Moreover, as in Corollary \ref{PdFamily}, it follows that $\widetilde \cU'$ is a coherent sheaf, and 
that {\bf if} $\widetilde \cU' = L_{\sigma'}^*\cU'$ for some object $\cU' \in \cD(S\times \cM')$, then 
$\cU'$ is a coherent sheaf, as well. The fact that $\cU'$ then parametrizes all $t_d-\epsilon$ stable 
objects was established in the interlude above.
But, as explained to us by Andrei C\u ald\u araru \cite{Cald}, the descent is an immediate consequence of Orlov's orthogonal decomposition of the derived category of the blow-up $\sigma'$ \cite{Orl92}. Indeed, it follows from the decomposition that {\it any} object $\cE \in \cD(S\times \widetilde \cM)$ whose restriction to the 
exceptional divisor descends in the derived category (i.e. $L{i_{D_S}}^*\cE \cong \sigma'^*\cF$ for some 
$\cF \in \cD(S\times \bP^\vee_d)$) must itself descend. 

\qed

\medskip

{\bf Proof of Step 3:} To recap, the moduli functor:
$$\{\mbox{flat families of $t_d-\epsilon$-stable objects of $\cA^\#$ with invariants
$(0,H,H^2/2)$}\}/\mbox{iso}$$
has the following properties:

\medskip

$\bullet$ the stable objects of $\cA^\#$ have only the automorphisms $\bC^*\cdot \mbox{id}$.

\medskip

$\bullet$ the proper variety $\cM'$ is in bijection with the set of stable objects, and

\medskip

$\bullet$ the coherent sheaf $\cU'$ on $S\times \cM'$ realizes this bijection.

\medskip

We want to conclude that $\cM'$ is a (fine) moduli space representing the functor, and 
$\cU'$ is a (universal) Poincar\'e object, which is well-defined up to tensoring by 
a line bundle pulled back from $\cM'$. This will follow
once we establish the a priori existence of an Artin {\it stack} for the functor: 
flat families of $t$-stable objects of $\cA^\#$ (with the given invariants, for each $t > \frac 16$, $t \ne t_d$).

\medskip

First, we appeal to Example (1.2) of the Appendix to conclude that functor: flat families of objects of
$\cA^\#$ with invariants $(0,H,H^2/2)$ is an Artin stack, which we will denote by
$\mathfrak M_{\cA^\#}(0,H,H^2/2)$. 
Step 3 is complete once we show that $t$-stability is an open 
condition on this functor, hence itself represented by an open substack. To this end, we prove first that 
$t$-stability for {\it some} $t > \frac 16$ is an open condition. 

\medskip

Suppose $B$ is a quasi-projective base scheme, and 
$\cE \in \cD(S\times B)$ is a flat family of objects of 
$\cA^\#$ with the given invariants. We may assume without loss of generality (see the proof of 
Lemma \ref{Sheaf} above) that $\cE$ is represented by a two-term complex $f:K \fun V$ where
$K,V$ are coherent sheaves on $S\times B$ such that $V$ is locally free and $K$ is flat (as a coherent 
sheaf) over $B$. By Theorem \ref{tstable}, if it is the case that $\cE_{S\times \{b\}}$ is $t$-stable 
for some closed point $b\in B$ and $t > \frac 16$, then either:
$f_b := f|_{S\times \{b\}}:K_b \fun V_b$ is injective and coker$(f_b) \cong i_*L_C$, or else:
$E_b = [K_b \stackrel f\fun V_b]$ fits in a short exact sequence (in $\cA^\#$) of the form 
$0 \fun \cI_W^\vee[1] \fun E_b \fun \cI_Z(H) \fun 0$. In the latter case, it follows in particular 
that $f_b$ fits in a long exact sequence:
$$0 \fun \cO_S \fun K_b \stackrel {f_b}\fun V_b \fun Q \fun 0$$
where $Q$ is a coherent sheaf fitting into: $0 \fun \cO_W^\vee[2] \fun Q \fun \cI_Z(H) \fun 0$.

\medskip

A little analysis gives the following {\it  necessary} conditions for $E_b$ to be $t$-stable for some $t > 1/6$:

\medskip

(a) $K_b$ is locally free (and then it follows that $\ker(f_b)$ is locally free).

\medskip

(b) $\ker(f_b) = H^{-1}(E_b)$ has rank $\le 1$, and $c_1(\ker(f_b)) \le 0$.

\medskip

(c) $len(tors(\coker(f_b))) < \frac {H^2}8$.

\medskip

On the other hand, (a)-(c) are very nearly {\it sufficient} conditions for $E_b$ to be $t$-stable for some $t > 1/6$. 
First,
if $f_b$ is injective, then clearly $E_b \cong i_*L_C$ where $L_C$ is torsion-free. If, on the other hand, 
$\ker(f_b)$ is a line bundle, then by (b), it must be of the form $\cO_S(-nH)$ for some $n \ge 0$. 
But $E_b \in \cA^\#$ has invariants $(0,H,H^2/2)$, by assumption, and this implies that $n = 0$, i.e. 
$\ker(f_b) \cong \cO_S$ and $\coker(f_b)$ has torsion only in dimension zero. Moreover, $\coker(f_b)/tors \cong \cI_Z(H)$ for some $Z \subset S$ of 
length equal to the length of the torsion (hence less than $H^2/8$), and as in Lemma \ref{Dual}, one may conclude that the kernel 
(in $\cA^\#$) of the map $E_b \fun \coker(f_b)/tors$ is of the form $\cI_W^\vee[1]$. Such extensions are necessarily $t$-stable for $t$ very close to the value $t_d$, where $d = len(Z)$, {\it provided that they are non-split}. 
Thus to ensure stability for some $t > 1/6$, we need only add:

\medskip

(d)  $E_b$ is not in the image of any of the (proper) morphisms: $h_d: S[d] \times S[d] \fun \mathfrak{M_{\cA^\#}}(0,H,H^2/2)$; $(W,Z) \mapsto \cI_W^\vee[1] \oplus \cI_Z(H)$ for any $d < H^2/8$.

\medskip

Finally, suppose $t_{d_0} < t < t_{d_0-1}$. Then $E_b$ is $t$-stable if it is $t$-stable for {\it some} $t > 1/6$, and 
moreover, it avoids both the images of the (proper) morphisms from Proposition \ref{Pd}:
$$i_d: \bP_d \fun \mathfrak{M_{\cA^\#}}(0,H,H^2/2) \ \mbox{defined by}\ \cE_d \ \ \mbox{for} \ \ d < d_0$$
and
$$i'_d:\bP_d^\vee \fun \mathfrak{M_{\cA^\#}}(0,H,H^2/2) \ \mbox{defined by}\ \cF_d \ \ \mbox{for} \ \ d \ge d_0$$

This completes the proof of the Step 3, and hence of the Theorem \qed.

\medskip

\rem It certainly ought to be possible to extend this theorem to remove the $t > 1/6$ assumption, i.e.
to produce Mukai flops for ``higher rank'' walls, in addition to the rank one walls. There are two 
places where improvements would need to be made to the proof. First, our argument in Step 3 for 
the openness of stability breaks down when higher rank walls are crossed (but a very recent result of 
Toda \cite{Tod07} gives an independent proof of the existence of the Artin stack). 
Second, there are cases where higher rank walls
coincide (with each other or with a rank one wall). In that case, the wall-crossing will not be a simple
Mukai flop, but more likely a stratified elementary modification of the sort investigated by Markman \cite{Mar01}
for the birational transformations of moduli spaces induced by Fourier-Mukai transforms.
In any case, this theorem is quite likely to generalize in many interesting ways.

\section{K3 Surfaces}

Let $S$ be a K3 surface of genus $g$ with Pic$(S) = \bZ[H]$ (i.e. $H^2 = 2g-2$). The results of 
\S 4 and \S 5 then give the following:

\medskip

\nt {\bf Vanishing:} For all $d < \frac {H^2}8 = \frac 14{(g-1)}$ and all pairs $Z,W \subset S$ of 
closed subschemes of length $d$,
$$H^i(S,\cI_W\otimes \cI_Z(H)) = 0 \ \mbox{for all $i > 0$}$$
(Proposition \ref{Van})

\nt {\bf Moduli Spaces:} By Theorem \ref{moduli}, there are Mukai flops of the relative Jacobian:
$$\cM_S\left(0,H,\frac{H^2}2\right) = \cM_0 --> \cM_1 --> \cdots --> \cM_{d_g}$$
as one crosses walls $t_d := \sqrt{\frac 14 - d\frac 2{H^2}} > \frac 16$. Thus, the final flop predicted by 
the Theorem finishes  with $\cM_{d_g}$ where $d_g = \lceil \frac{H^2}9\rceil$ is the round-up of 
$H^2/9 = (2g-2)/9$.

\medskip

However, these results for K3 surfaces are not optimal. For example, the vanishing theorem
predicts that $\cO_S(H)$ is very ample on a K3 surface of genus $\ge 6$, whereas in fact 
it is very ample for genus $\ge 3$. We can get better results if we use Bridgeland's central charge:
$$Z'_{(D,F)}(E) := -\int_S e^{-(D + iF)} \ch([E]) \sqrt {\mbox{td}(S)} = Z_{(D,F)}(E) - \rk(E)$$
instead of the central charge (without the todd class contribution) $Z_{(D,F)}$ of \S 2. The key point is that if $E$ is an
$H$-stable vector bundle on a K3 surface $S$, then:
$$\chi(S,E\otimes E^*) \le 2$$
and a quick computation with Riemann-Roch then shows that: 
$$\ch_2(E) \le \frac{c_1^2(E)}{2\cdot \rk(E)} - \rk(E) + \frac 1{\rk(E)}$$
(for all $H$-stable torsion-free sheaves), which is sharper than the Bogomolov-Gieseker
inequality. For the choice $Z'_t :=  Z'_{(\frac 12 H, tH)}$ and category $\cA^\#$ as before, we now get the following version of 
Corollary \ref{slope}:

\begin{cor}\label{slopeprime} On a K3 surface $S$ of genus $g$ and Picard number one, $(Z'_t,\cA^\#)$ is a Bridgeland slope function if either
$g$ is odd and $t > 0$, or else $g$ is even and $t > \sqrt{1/(4g-4)}$.
\end{cor}

\pf The proof of Corollary \ref{slope} holds up until the final case:
$$\mu_{tH}(E) = \left(\frac 12 H\right)\cdot (tH) \Leftrightarrow \mu_H(E) = 
\frac{c_1(E)\cdot H}{\rk(E)} = \frac 12 H^2$$
which, because of the Picard number one assumption, implies that:
$$c_1(E) = cH, \ \ \ \ \rk(E) = 2c \ \mbox{for some integer $c$}$$
For torsion-free sheaves with these invariants, we then use the improved inequality to obtain:
$$\mbox{Re}(Z'_t(E)) \ge t^2 \rk(E)\frac{H^2}2 - \frac{1}{\rk(E)} = t^2cH^2 - \frac 1{2c}$$
and this is positive if $t^2(H^2) > 1/2$ or $t > \sqrt{1/(4g-4)}$, giving the result for even genera.

\medskip

When $g$ is odd and $r = 2c$ is even, then the inequality can be improved:
$$\ch_2(E) \le \frac{g-1}{2c} - 2c + \frac 1{2c} \Rightarrow \deg(\ch_2(E)) \le \frac{g-1}{2c} - 2c$$
because the $\deg(\ch_2(E))$ is an integer(!) This gives the better result for odd genera.

\qed

\nt Lemma \ref{LambdaI} remains valid (with the same proof!) for $Z'_t$ and we obtain the following:

\medskip

\begin{pro} [\bf Better Vanishing]\label{Vanprime} For a K3 surface $S$ of genus $g$ with 
$\emph{Pic}(S) = \bZ[H]$, 
$$H^i(S,\cI_W\otimes \cI_Z(H)) = 0 \ \ \mbox{for all $i > 0$ and all subschemes $Z,W$ of length $d
< \frac{g+2}4$}$$
\end{pro}

\pf As in the proof of Proposition \ref{Van}, it suffices to exhibit a value of $t$ such that 
$(Z'_t,\cA^\#)$ is a Bridgeland  slope function (i.e. satisfies the criteria of Corollary \ref{slopeprime}) 
and:
$$\mu'_t(\cI_W^\vee[1]) = 0 =  \mu'_t(\cI_Z(H))$$
where $\mu' = -\mbox{Re}(Z')/\Im(Z')$. Solving for this equality yields:
$$d = \frac {H^2}8 - \frac{t^2H^2}2 + 1 = \frac{g+3}4 - t^2(g-1)$$
and keeping mind the constraints on $t$, this yields $d < \frac{g+3}4$ if $g$ is odd, and 
$d < \frac{g+2}4$ if $g$ is even. Since $d$ must be an integer, the result of the Proposition is 
then sharp.

\qed

\begin{rem} This bound is now the best possible, since, for example, it gives very ampleness 
of $\cO_S(H)$ in 
genus $3$, but just fails to give it in genus $2$ (where $\cO_S(H)$ is NOT very ample!).

\end{rem}

To get the right number of flops on a K3 surface, note that the flops will exist for all the critical values:
$$t = \sqrt{\frac{\frac{g+3}4 - d}{g-1}}$$
until $t$ runs into the first ``higher (odd) rank wall,'' (the value $t = \frac 16$ in the proof of Theorem 
\ref{tstable}). As in that Theorem, this highest value of $t$ is computed by setting:
$$\mbox{Re}(Z'_t(E)) = 0$$
where $E$ is an $H$-stable torsion-free sheaf where $\rk(E) = 3$, $c_1(E) = 2H$, and $\ch_2(E)$ is maximal. Using the inequality:
$$\deg(\ch_2(E)) \le \frac{(2H)^2}{6} - 3 + \frac 13 = \frac {4g}3 - 4$$
and the fact that it is an integer, we obtain (sharp) values for such $t$ by setting:
$$0 = \left\{
\begin{array}{ll} -\frac 1{24}H^2 + \frac 32 t^2H^2 - \frac 13 & \mbox{if}\ \ g = 3n \\ \\
-\frac 1{24}H^2 + \frac 32 t^2H^2  & \mbox{if}\ \ g = 3n + 1 \\ \\
-\frac 1{24}H^2 + \frac 32 t^2H^2  + \frac 13& \mbox{if}\ \ g = 3n + 2 
\end{array}\right.$$

It follows that in all genera $> 2$, the upper bound on the number of Mukai flops improves to:
$$d_g = \lceil \frac {2(g+3)}9\rceil$$

\medskip

\nt {\bf Small Genus Examples:} {\bf Genus 2.} In all genera $\ge 2$, the first Mukai flop exists:
$$\cM_0 --> \cM_1$$
replacing
$\bP^g := \bP(H^0(S,\cO_S(H)))$ with its dual $(\bP^g)^\vee = \bP(\Ext^2_S(\cO_S(H),\cO_S))$. 

\medskip

This Mukai flop can also be realized with the (standard) Fourier-Mukai transform:
$$i_*L_C \mapsto R{\pi_2}_*(R\pi_1^*i_*L_C \stackrel L \otimes \cI_\Delta)^\vee$$
where $\pi:S\times S \fun S$ are the projections, and $\Delta \subset S\times S$ is the diagonal.

\medskip

We claim that, in fact, 
$\cM_1 \cong \cM_S(g-1,H,-H^2/2)$. 
Here is a sketch of the argument, which is well-known \cite{Mar01}. Consider first the case $L_C \not\cong \omega_C$. Then $\dim(H^0(S,i_*L_C)) = g-1$, and:
$$E = [H^0(S,i_*L_C)\otimes \cO_S \stackrel r \fun i_*L_C]^\vee$$
is a stable torsion-free sheaf with the desired invariants (which is locally free iff  $r$ is surjective). This gives a rational 
map $\Phi: \cM_0 --> \cM_S(g-1,H,-H^2/2)$ that is well-defined off the locus $\bP^g$
parametrizing the sheaves of the form $i_*\omega_C$. To define $\Phi$
further, we blow up $\bP^g \subset \cM_0$ by choosing codimension one subspaces $V \subset H^0(C,\omega_C)$
at each point $i_*\omega_C \in \bP^g$. Each such point is then mapped to $[V\otimes \cO_S \fun i_*\omega_C]^\vee$, extending $\Phi$ to $\widetilde \cM$, and furthermore, one can check that the 
map $\Phi$ blows down the exceptional divisor $D$, giving the Mukai flop. Notice that in the 
genus $2$ case, $\cM_S(1,H,-H^2/2) \cong S[2]$ is the Hilbert scheme parametrizing sheaves of 
the form $\cI_Z(H)$.

\medskip

{\bf Genus 3.} In all genera $\ge 3$, there is a second Mukai flop:
$$\cM_1 --> \cM_2$$
In genus $3$, the flop locus $\bP_1 \subset \cM_1$ is a {\it divisor} (codimension $3-2(1) = 1$).
Thus in genus $3$, this second ``flop,''  which occurs at $t = \frac 12$ (corresponding to $d=1$) is actually the identity. On the other hand, there is one further ``rank three wall'' at $t = \sqrt \frac 1{12}$, and 
the indeterminacy locus for that wall is, again, $(\bP^g)^\vee$ so that the rank three wall 
 is the inverse of the original flop(!) That is, there is a symmetry:
$$\cM_0  \stackrel{\rm 1st\ rank\ one\ flop}{-->} \cM_1 \stackrel{\rm 2nd\ rank\ one\ flop}{\cong} \cM_1 
\stackrel{\rm rank\ 3\ flop}{-->} \cM_0$$

Experimental evidence suggests some sort of similarly symmetric picture when {\it all} the 
higher rank walls are taken into account in each genus  of the form $4n + 3$.

\medskip

\section{Abelian Surfaces}

Let $S$ be a simple abelian surface, with $(1,D)$ polarization and NS$(S) \cong Z[H]$, so:
$$H^2 = 2D \ \ \mbox{and}\ \ h^0(S,\cO_S(H)) = D$$

Then the vanishing and main theorem hold in the following modified forms:

\medskip

\nt {\bf Vanishing:} For all subschemes $Z,W \subset S$ of length $< H^2/8 = D/4$ and all $\cL \in \Pic^0(S)$, 
$$\mbox{H}^i(S,\cI_W\otimes \cI_Z(H) \otimes \cL) = 0 \ \mbox{for}\ \ i = 1,2$$

No better vanishing is to be expected, since the Bogomolov-Gieseker bound matches the bound coming from
Riemann-Roch applied to:
$$\chi(E\otimes E^*) \le 2 \ \mbox{for $H$-stable bundles}$$ 
In this case, as well, td$(S) = 1$, and the two central charges $Z$ and  $Z'$ coincide.

\medskip

\nt {\bf Moduli:} Let $\widehat S = \mbox{Pic}^0(S)$ be the dual abelian surface. 
Then Proposition \ref{Pd} here can be easily modified to produce projective bundles:
$$\bP_d, \bP_d^\vee \ \mbox{over}\ \ (\widehat S\times S[d]) \times (\widehat S\times S[d])$$
parametrizing extensions:
$$0 \fun \cL_1\otimes \cI_Z(H) \fun E_d \fun \cL_2\otimes \cI_W^\vee[1] \fun 0 \ \mbox{and}\ 
0 \fun \cL_2\otimes \cI_W^\vee[1] \fun F_d \fun \cL_1\otimes \cI_Z(H) \fun 0$$
(for $\cL_1,\cL_2 \in \widehat S$),
which embed:
$$\bP_d \hookrightarrow \cM_{t_d + \epsilon} \  \ \ \mbox{and}\ \ \ \bP_d^\vee \hookrightarrow \cM_{t_d - \epsilon}$$
as before, as the centers for the birational transformations.

\medskip

Here, the dimension count is:
$$\dim(\cM) = 2 - \chi(\mbox{RHom}_{\cD(S)}(E,E)) = 2 + H^2 = 2D + 2$$
and $\bP_d$ has fiber dimension $D - 2d - 1$ over a base of dimension $4d+4$, which therefore continues to satisfy:
$$\mbox{(fiber dimension)} = \mbox{(codimension)},$$
the condition for the birational transformation replacing $\bP_d$ with $\bP_d^\vee$ to be a Mukai flop.

\medskip

And the number of such flops is $\lceil \frac{H^2}9\rceil = \lceil \frac {2D}9 \rceil$, as before, as well.

\medskip

\nt {\bf Small Values for D:} 

\medskip

{\bf D = 1,2,3,4.} Only the $d = 0$ flop exists. In the $D = 1$ (principal polarization) 
case, the map:
$$\bP_0 = \widehat S\times \widehat S \fun \cM$$
is an isomorphism. 

\medskip

{\bf D = 5.} Vanishing for $d = 1$ gives the sharp (and well-known) result:
$$H \ \mbox{(and its translates) are very ample when} \ D \ge 5$$
In this $D = 5$ case, the embedding $\bP_1 \hookrightarrow \cM \ \mbox{has codimension two}$.

\section{Stable Pairs}

Recall that for a smooth projective curve $C$ with (ample) line bundle $L$, Serre duality gives:
$$\bP := \bP(H^0(C,\omega_C\otimes L)^*) \cong \bP(\Ext^1_C(L,\cO_C))$$
which is then on the one hand a projective space for extensions:
$$0 \fun \cO_C \fun E \fun L \fun 0$$
and on the other, the image for the ``linear series map''  $\phi_{L\otimes \omega_C} : C \fun \bP$.

\medskip

Michael Thaddeus \cite{Tha94} showed that there is a one-parameter family of stability conditions on isomorphism 
classes of ``pairs''
$(\cO_C \fun E)$ satisfying $\det(E) \cong L$ that exhibit wall-crossing behavior and a sequence of birational moduli spaces:
$$\bP =: \cP_0 --> \cP_1 --> \cP_2 --> \cdots --> \cP_{\lfloor \frac{d-1}2 \rfloor}; \ \ d = \deg(L)$$
such that:

\medskip

$\bullet$ $\cP_1$ is the blow-up of $\cP_0$ along the embedded curve $C$.

\medskip

$\bullet$ $\cP_{d+1}$ is a ``flip'' of $\cP_d$, replacing the proper transform of the secant variety of 
projective $d$-planes spanned by $d+1$ points of $C$
with a projective bundle:
$P^\vee_d \fun \mbox{Sym}^{d+1}(C)$ with fiber
$$\bP(H^0(C,\omega_C\otimes L(-2D))^*) = \bP(\Ext^1_C(L(-D),\cO_C(D))) \ \mbox{over}\  D \in \mbox{Sym}^{d+1}(C)$$

\medskip

In \cite{Ber97}, the second author asked whether such a sequence of 
moduli spaces might also exist for embeddings
$\phi_{L\otimes \omega_C}:X \rightarrow \bP$ by a sufficiently ample line bundle $L$. We remark here that 
a non-trivial sequence of ``Thaddeus-like'' flips does indeed exist when $S$ is a simple K3 surface of genus $g$, by restricting the chain of Mukai flops above:
$$\cM_1 --> \cM_2 --> \cdots --> \cM_{\lceil \frac {2g+6}9\rceil}$$
to the subvariety:
$$\bP_0^\vee = \bP(H^0(S,\cO_S(H))^*) = \bP(\Ext_S^2(\cO_S(H),\cO_S)) \subset \cM_1$$

\medskip

\nt {\bf Definition:} We will call the proper transform,  $\cP_d \subset \cM_{d+1}$, of $\bP_0^\vee \subset \cM_1$, 
the space of {\it stable pairs} for $t$-stable objects ($t_d < t < t_{d+1}$) with invariants $(0,H,H^2/2)$.

\medskip

We now describe candidates for the points of $\cP_d \subset \cM_{d+1}$ (without proof) in a rather close analogy
with the  case of Thaddeus stable pairs (on curves).

\medskip

$\bullet$ Each point of $\cP_d$ represents an object $E \in \cA^\#$ with a (unique!) non-trivial ``section:''
$$\cO_S[1] \rightarrow E$$

Indeed, {\bf every} one of the objects parametrized by each $\bP_d^\vee$ fits in an exact sequence:
$$0 \fun \cI_W^\vee[1] \fun E \fun \cI_Z(H) \fun 0$$
and thus satisfies $H^{-1}(E) = \cO_S$, so that the canonical inclusion $H^{-1}(E)[1] \fun E$ is a ``section.''
Moreover, since Hom$(\cO_S[1], H^0(E)) = \Ext_S^{-1}(\cO_S,H^0(E)) = 0$, it follows that this section is unique (up to scalars). Thus the space 
of $t$-stable objects $E\in \cM_d$ admitting sections is the union of the proper transforms of 
all the $\bP^\vee_e \subset \cM_{e+1}$ for all $e < d$. In particular, it contains $\cP_d$ as a component. 

\medskip

$\bullet$ The quotient by the section, $E/\cO_S[1]$, satisfies 
$$0 \fun \cO_Z^\vee[2] \fun E/\cO_S[1] \fun \cI_Z(H) \fun 0 \ \mbox{for some $Z \subset S$ of length $\le d$}$$
(Note: $\cO_Z^\vee[2]$ is the torsion coherent sheaf ${\ms E}xt_S^2(\cO_Z,\cO_S)$ of length
$len(Z)$ supported on $Z \subset S$).
This is analogous to fixing the determinant $L = \cO_C(H)$ and noting that the section $\cO_C \rightarrow E$, if it vanishes
along an effective divisor $D \subset C$, gives rise to an exact sequence for the quotient:
$$0 \fun \cO_D^\vee[1] \fun E/\cO_C \fun L_C(-D) \fun 0$$
where $\cO_D^\vee[1]$ is the torsion coherent sheaf ${\ms E}xt^1_S(\cO_D,\cO_C)$ 
of length $\deg(D)$ supported on $D$.

\medskip

$\bullet$ The sequence above splits.

\medskip

This last condition is automatic for curves, by the classification of modules over a PID, but not so in the surface case. Indeed, fixing $Z \subset S$ of length $d+1$, the space of extensions of the form:
$$0 \fun \cI_Z^\vee[1] \fun E \fun \cI_Z(H) \fun 0$$
is parametrized by $\bP(H^0(S,\cI_Z\otimes \cI_Z(H))^*)$, and one can check that the splitting of the sequence for 
$H^0(E)/\cO_S[1]$ restricts the extensions to be of the form:
$$H^0(S,\cI_Z^2(H))^* \subset H^0(S,\cI_Z\otimes \cI_Z (H))^*$$
These are the extensions that would be ``inserted'' at the $d$th (Thaddeus) flip.

\medskip

This is now in perfect analogy with the projective bundle $P^\vee_d \fun \mbox{Sym}^{d+1}(C)$
that appears in the Thaddeus flips for curves, though 
in this case it will not be a projective bundle, since the dimensions of the spaces $H^0(S,\cI_Z^2(H))$ 
will jump up.

\medskip

\nt {\it Final Remark:} We hope to be able to  ``fix'' this definition by constructing the 
appropriate moduli problem for stable pairs on a surface (not necessarily K3) in order to 
define $\cP_d$ as a moduli space, to determine its scheme structure through deformation theory.

\appendix

\section{The openness of tilted hearts (by Max Lieblich)}

\begin{abstract}
   We show that there is a natural condition on a torsion theory on the
   category of coherent sheaves on a flat proper morphism which ensures
   that the heart of the tilting is represented by an Artin stack
   locally of finite presentation over the base.
\end{abstract}

   Let $X\to S$ be a flat proper morphism of finite presentation
   between schemes.  Write $\cA\to S$ for the fpqc stack of
   quasi-coherent sheaves on $X$, $\cA_{pf}$ for the substack
   parametrizing quasi-coherent sheaves of finite presentation, and
   $\cA_{pf}^p$ for the substack of $\cA_{pf}$ consisting of
   families of quasi-coherent sheaves of finite presentation which are
   flat over the base.  For the moment, we work with the full stacks of
   categories and not merely the underlying stacks of groupoids.  It is
   a standard result that the stack of groupoids
   $(\cA_{pf}^p)^{\text{\rm gr}}$ underlying $\cA_{pf}^p$ is an Artin
   stack locally of finite presentation over $S$.

\begin{rem}
   Note that while $\cA$ is a stack of abelian categories, this is
   not in fact true of $\cA_{pf}$.  A simple example which shows that
   $\cA_{pf}$ is not abelian is the following: the homomorphism of
   finitely presented $\bZ[x_1,x_2,\dots]$-modules
   $\bZ[x_1,x_2,\dots]\to\bZ[x_1,x_2,\dots]$ which sends $x_{2i-1}\mapsto
   x_{2i-1}$ and $x_{2i}\mapsto x_{2i-1}$ (for $i\geq 1$) has kernel
   $(x_{2i}-x_{2i-1})$, which is not even finitely generated.
   Moreover, this difficulty cannot be avoided by requiring $S$ to be
   Noetherian, since in the theory of stacks one must allow arbitrary base
   changes.

   However, for any field-valued point $s\to S$, the fiber category
   $(\cA_{pf})_{s}$ is abelian: it is the category of coherent
   sheaves on a finite type $\kappa(s)$-scheme.
\end{rem}

\begin{defi}
   A \emph{stack of torsion theories\/} in $\cA_{pf}$ consists of a
   pair of full substacks $(T,F)$ of $\cA_{pf}^p$ with the property
   that for each point $s=\Spec K\to S$, the pair of subcategories
   $(T_{s},F_{s})$ in $(\cA_{pf})_{s}$ is a torsion theory in
   the classical sense.  A stack of torsion theories $(T,F)$ is
   \emph{open} if the groupoids underlying $T$ and $F$ are open
   substacks of $(\cA_{pf}^p)^{\text{\rm gr}}$.
\end{defi}

\begin{lem}
   Let $(T,F)$ be a stack of torsion theories.  Suppose $\Spec
   K\to S$ is a point and $L/K$ is a field extension.  An object
   $M\in\cA_K$ is in $T_K$ if and only if $M|_L$ is in $T_L$, and
   similarly for $F$.
\end{lem}
\pf Since $T\subset\cA$ is a full fpqc substack and $L/K$ is
   faithfully flat, the result is immediate: an object $t\in T_L$
   acquires a descent datum relative to $K$ by transport of structure
   from $M_L$.  This of course works similarly for $F$.
\qed

\medskip

\noindent Thus, belonging to $T$ or $F$ is determined by geometric fibers.

\medskip

Lest the reader give up reading in disgust, let us give a couple of
examples when $X$ is a projective variety over an algebraically closed
field.

\begin{example}\label{ex:conds}
The two most important examples (for our purposes) are
   the following.
   \begin{enumerate}
   \item If $X$ is a smooth projective variety, then letting $T$ be the
     stack of torsion sheaves and letting $F$ be the stack torsion free
     sheaves on $X$ (i.e., pure sheaves on $X$ of maximal dimension)
     defines an open stack of torsion theories.
   \item\label{item:1} If $X$ is a K3 surface, Bridgeland has described
     a class of examples (see Lemma 5.1 of \cite{bridgeland-k3}).  Given
     an ample divisor $\omega$ on $X$ and a class
     $\beta\in\NS (X)\otimes \bR$, let $T$ be the category
     of coherent sheaves on $X$ whose torsion free parts have the
     property that all subquotients of the Harder-Narasimhan filtration
     have slope strictly larger than $\beta\cdot\omega$, and let $F$ be
     the category of torsion free sheaves on $X$ whose
     Harder-Narasimhan factors all have slope at most
     $\beta\cdot\omega$.
   \end{enumerate}
   \pf[Proof that (\ref{ex:conds}.\ref{item:1}) defines open
     stacks of torsion theories]
     In the following, we will repeatedly use the fact that the torsion
     free locus of a flat family of quasi-coherent sheaves of finite
     presentation is open.  The proof is similar to those given here,
     and we leave it to the reader.  (It is conceptually somewhat
     easier to prove the equivalent assertion that the locus with
     non-trivial torsion subsheaf is closed.)  To show that the
     condition that every Harder-Narasimhan factor has slope at most
     $\beta\cdot\omega$ is open, we will use a standard argument: we
     will show that the locus is constructible and stable under
     generization.  More precisely, let $\ms F$ be a flat family of
     quasi-coherent sheaves of finite presentation on $X\times B$ with
     torsion free fibers.  We will show that the set of points
     $U\subset B$ over which the fibers of $\ms F$ are in $F$ is open.
     By standard limiting arguments, we may assume that $B$ is of
     finite type over $k$.

     To show that $U$ is constructible, we may assume that $B$ is
     reduced and irreducible, and we wish to show that $U$ contains the
     generic point if and only if it contains an open subset of points.
     This follows from the fact that the slope is constant in a flat
     family and the existence of the relative Harder-Narasimhan
     filtration over a dense open subscheme of $B$ (Theorem 2.3.2 of
     \cite{h-l}).

     To show that $U$ is stable under generization, let $R$ be a
     discrete valuation $k$-algebra and $\ms F$ a flat family of
     torsion free coherent sheaves on $X\otimes R$ such that the closed
     fiber $\ms F_s$ is in $F$.  The maximal destabilizing subsheaf
     $\ms G_{\eta}\subset\ms F_{\eta}$ on the generic fiber
     extends to a coherent subsheaf $\ms G\subset\ms F$ such that $\ms
     F/\ms G$ is $R$-flat.  It follows that the closed fiber $\ms G_s$
     gives a subsheaf of $\ms F_s$ whose slope must be at most
     $\beta\cdot\omega$.  Since the slope is constant in a flat family,
     we see that $\mu(\ms G_\eta)\leq\beta\cdot\omega$, as desired.
     (That this passes to the geometric generic fiber follows from the
     compatibility of the Harder-Narasimhan filtration of $\ms
     F_{\eta}$ with extension of the base field, Theorem 1.3.7 of
     \cite{h-l}.)

     The proof that $T$ is open is similar, but with an extra
     complication due to the presence of the torsion subsheaf.  The
     point is that in both the proof of constructibility and stability
     under generization, one can assume that the the torsion subsheaves
     form a flat subfamily of $\ms F$.  Thus, one immediately reduces
     to showing that for a torsion free family $\ms F$, the locus over
     which all Harder-Narasimhan factors have slope strictly greater
     than $\beta\cdot\omega$ is open.  One can again use the existence
     of the relative Harder-Narasimhan filtration to get
     constructibility.  Stability under generization can be proven by
     induction on the rank as follows.  Note that if the
     Harder-Narasimhan factors have slope larger than
     $\beta\cdot\omega$, then $\mu(\ms F_s)>\beta\cdot\omega$.  Thus,
     if $\ms F$ has semistable generic fiber then $\ms F_\eta$ must be
     in $T$ when $\ms F_s$ is.  On the other hand, if $\ms F_\eta$ is
     not semistable, then there is a flat subfamily $\ms G\subset\ms F$
     which agrees with the maximal destabilizing subsheaf on the
     generic fiber.  On the other hand, $\mu(\ms G_\eta)\geq\mu(\ms
     F_{\eta})=\mu(\ms F_s)>\beta\cdot\omega$.  The quotient $\ms
     F/\ms G$ is still flat, and the closed fiber is a quotient of $\ms
     F_s$.  Since $T_s$ is closed under the formation of quotients, we
     conclude by induction that $\ms F_\eta/\ms G_\eta$ is in $T_\eta$,
     whence $\ms F_\eta$ is in $T_\eta$.
   \qed
\end{example}

The formation of the derived category yields a fibered category
$\cD\to S$ which over $B\to S$ takes the value $D(\cA_B)$.  The fibered
category structure comes from the \emph{derived pullback\/} functors
(and their natural functorialities).  The substack $\cA_{pf}$ gives
rise to a substack $\cD_{pf}$ by taking $(\cD_{pf})_B$ to be the
subtriangulated category of $D(\cA_B)$ generated by complexes with
entries in $(\cA_{pf})_B$.  (Note that the subtriangulated category
generated by the same procedure by $\cA_{pf}^p$ is not the same,
although it does agree on fiber categories over field-valued points of
$S$.)


Recall that a complex $E$ on $X$ is \emph{relatively perfect\/} if for
every affine $\Spec A$ of $S$ and $\Spec B$ of $X$ such that $\Spec B$
maps into $\Spec A$ under the structural morphism $X\to S$, the
complex $E|_{\Spec B}$ is quasi-isomorphic to a bounded complex of
$A$-flat coherent $B$-modules.  The complex $E$ is \emph{universally
   gluable\/} if for every geometric point $\bar s\to S$ and every
$i>0$, $\Ext^i_{X_{\bar s}}(E_{\bar s},E_{\bar s})=0$.  Note that if
$E$ is in the heart of a sheaf of $t$-structures on $X$, then it must
be universally gluable.  The sub-fibered category of $\cD$ formed by
relatively perfect universally gluable complexes is denoted $\ms
D_{pug}(X/S)$.  It is a standard result (Corollaire 2.1.23 of
\cite{BBD} or Theorem 2.1.9 of \cite{abramovich-polishchuk}) that $\ms
D_{pug}(X/S)$ is a stack on the fpqc topology on the category of 
$S$-schemes.

We recall the main theorem of \cite{moduli-complexes}.

\begin{thm}
   The fibered category $\ms D_{pug}(X/S)\to S$ is an Artin stack
   locally of finite presentation.
\end{thm}

\noindent This is the Mother of All Moduli Spaces: it contains the
hearts of all of the sheaves of $t$-structures on $X$.  In our case,
we can make this precise as follows.

Given a stack of torsion theories $(T,F)$, we can define a substack
$\ms D_{(T,F)}(X/S)$ of $\ms D_{pug}(X/S)$ corresponding to the family of
hearts of the tilting with respect to the torsion theory.

\begin{defi}
   Given a stack of torsion theories $(T,F)$, the \emph{stack of tilted
     hearts with respect to $(T,F)$\/} is the stack $\ms
   D_{(T,F)}(X/S)$ whose objects over $B\to S$ are objects $\ms C$ of
   $\ms D_{pug}(X/S)_B$ such that for every geometric point $\bar s\to
   B$, the derived pullback $\ms C|_{\bar s}^{\mathbf L}\in\bD(X_{\bar
     s})$ is in the heart of the tilting with respect to the torsion
   theory $(T_s,F_s)$, i.e., it has cohomology only in degrees $-1$ and
   $0$ with $\ms H^{-1}(E)\in F_s$ and $\ms H^0(E)\in T_s$.
\end{defi}

The main result of this appendix is the following.

\begin{thm}
   If $(T,F)$ is an open stack of torsion theories then $\ms
   D_{(T,F)}(X/S)$ is an open substack of $\ms D_{pug}(X/S)$.
\end{thm}

\pf
   We will show something \emph{a priori\/} more general: given an
   affine scheme $B\to S$ and a relatively perfect complex $E$ on
   $X\times_S B$, there is an open subscheme $U\subset B$ parametrizing
   fibers in the heart of the tilting with respect to $(T,F)$.  More
   precisely, we will show that there is an open subscheme $U$ such
   that for a point $b\to B$, the derived base change $E_{b}$ is in
   $\ms D_{(T,F)}$ if and only if $b$ factors through $U$.  It follows
   from 2.2.1 of \cite{moduli-complexes} and 8.10.5 of \cite{ega4-3}
   that we may assume $B$ is Noetherian.

   Let $U\subset B$ be the subset parametrizing points over which the
   fiber of $E$ is in $\ms D_{(T,F)}$.  To show that $U$ is open, it
   suffices to show that $U$ is constructible and stable under
   generization.  By standard results (e.g., 2.1.3 and 2.1.4 of
   \cite{moduli-complexes}), the locus in $B$ over which the cohomology
   of the geometric fibers of $E$ is concentrated in degrees $-1$ and
   $0$ is open.  Thus, we may assume from the start that $\ms H^i(E)=0$
   unless $i\in\{-1,0\}$.

   Since $B$ is Noetherian, it follows from the results of \S 9.2 of
   \cite{ega3-1} that to show $U$ is constructible, it suffices to
   assume that $B$ is reduced and irreducible, and then we wish to show
   that the generic point of $B$ is in $U$ if and only if an open
   subset of points is contained in $U$.  We may thus shrink $B$ until
   the (coherent) cohomology sheaves $\ms H^0(E)$ and $\ms H^{-1}(E)$
   are $B$-flat.  It now follows from the standard spectral sequences
   that the formation of $\ms H^0$ and $\ms H^{-1}$ are compatible with
   arbitrary base change on $B$.  Since $T$ and $F$ are open, we see
   that $U$ must be open.

To show that $U$ is stable under generization, we may assume that
$B=\Spec R$ is the spectrum of a discrete valution ring and that the
special fiber is in $\ms D_{(T,F)}$.  Write $b$ for the closed point
of $B$ and $\eta$ for the open point.  Let $t\in R$ be uniformizer.
There is an exact sequence
$$0\to\ms H^{-1}(E)\xrightarrow{t}\ms H^{-1}(E)\to
   \ms H^{-1}(E_b)\to\ms H^0(E)\xrightarrow{t}\ms H^0(E)\to\ms
   H^0(E_b)\to 0,$$
where the indicated arrows are multiplication
by $t$.  We see immediately that $\ms H^0(E_b)=\ms
H^0(E)\otimes\kappa(\eta)$; moreover, $T$ is closed under the formation
of quotients, so any quotient of $\ms H^0(E)\otimes\kappa(\eta)$ lies in
$T_b$.  Dividing out $\ms H^0(E)$ by its associated submodules lying
over $b$ yields an $R$-flat quotient sheaf $\ms Q$ with generic fiber
$\ms H^0(E_\eta)$ and such that $\ms Q_b$ is a quotient of $\ms
H^0(E_b)$.  Thus, $\ms Q_b$ lies in $T_b$, whence the generic fiber
$\ms H^0(E_\eta)$ must lie in $T_\eta$, as $T$ is an open substack of
the stack of \emph{flat\/} families of coherent sheaves.

To prove that $\ms H^{-1}(E_\eta)\in F_\eta$ is somewhat simpler.  It
follows from the exact sequence that $\ms H^{-1}(E)$ is flat over $R$;
since $F$ is closed under subobjects and $\ms H^{-1}(E_b)\in F_b$, the
openness of $F$ (which is a substack of the stack of flat
families) shows that $\ms H^{-1}(E)\in F_B$.  Thus, $\ms
H^{-1}(E_\eta)\in F_\eta$, as required.
\qed

\begin{cor}
   The fibered category $\ms D_{(T,F)}(X/S)\to S$ is an Artin stack locally
   of finite presentation.
\end{cor}

\end{document}